\documentclass[12pt]{article}
\usepackage{amssymb}
\usepackage{amsmath}
\usepackage{latexsym}
\usepackage{mathrsfs}

\numberwithin{equation}{section}
\newtheorem{thm}[equation]{Theorem}

\newtheorem{defn}[equation]{Definition}
\newtheorem{rem}[equation]{Remark}
\newtheorem{lem}[equation]{Lemma}
\newtheorem{corol}[equation]{Corollary}

\title{Sharp spectral stability estimates via the Lebesgue measure of domains for higher order elliptic operators\footnote{to appear in {\it Revista Matem\'{a}tica Complutense}, DOI 10.1007/s13163-011-0079-2}}
\author{Victor I.  Burenkov and Pier Domenico Lamberti\footnote{Corresponding Author}}

\date{ \today }

\begin{document}

\newcommand{\rea}{\mathbb{R}}

\maketitle

%
%
%

\noindent
{\bf Abstract:} We prove sharp stability estimates for the variation of the eigenvalues of non-negative self-adjoint elliptic operators
of arbitrary even order upon variation of the open sets on which they are defined. These estimates are expressed in terms of the Lebesgue measure of the symmetric difference of the
open sets. Both Dirichlet and Neumann boundary conditions are considered.

\vspace{11pt}

\noindent
{\bf Keywords:} Elliptic equations, Dirichlet and Neumann
boundary conditions,
stability of eigenvalues, sharp estimates,
domain perturbation.

\vspace{6pt}
\noindent
{\bf 2000 Mathematics Subject Classification:} 35P15, 35J40, 47A75, 47B25.
%
%



\section{Introduction}

We consider a non-negative self-adjoint operator

\begin{equation}
\label{classic}
Hu=
 (-1)^m\sum_{|\alpha |=|\beta |=m} D^{\alpha }\left(A_{\alpha \beta}(x)D^{\beta }u  \right),\
\ \ x\in\Omega ,
\end{equation}
of order $2m$ subject to homogeneous Dirichlet or Neumann  boundary conditions on a bounded  open set $\Omega $
in ${\mathbb{R}}^N$.  Here $m\in {\mathbb{N}}$ is arbitrary and  the coefficients $A_{\alpha \beta }$ are bounded measurable  functions satisfying the uniform ellipticity condition.

If $\Omega $ is sufficiently regular then  $H$ has compact resolvent and its spectrum consists of a sequence of eigenvalues
$$
\lambda_1[\Omega ]\le \lambda_2[\Omega ]\le \dots \le \lambda_n[\Omega]\le \dots
$$
of finite multiplicity such that $\lim_{n\to \infty }\lambda_n[\Omega ]=\infty $.

In this paper, for fixed coefficients $A_{\alpha \beta}$, we prove sharp stability estimates for the variation of $\lambda_n[\Omega ]$ upon variation
of $\Omega$.

The problem of estimating the deviation of the eigenvalues of second order elliptic operators following a domain perturbation  has been considered by several  authors: we refer to Burenkov, Lamberti and Lanza de Cristoforis~\cite{burlamlan} for extensive references on this subject and to Barbatis, Burenkov and Lamberti~\cite{babula} for a recent paper
concerning stability estimates for resolvents, eigenfunctions and eigenvalues in the  case of domain perturbations obtained by suitable diffeomorphisms.

The case of higher order operators has been far less investigated. We refer to Prikazhchikov and Klunnik~\cite{priklu} for the case of the biharmonic operator subject to Dirichlet boundary conditions on smooth open sets and to Burenkov and Lamberti~\cite{bulahigh} for the general case of higher order elliptic operators  subject to Dirichlet or Neumann boundary conditions on open sets with continuous boundaries. The estimates provided in \cite{bulahigh},\cite{priklu} are expressed in terms
of the Hausdorff distance between the open sets.

In this paper we develop  the approach of Burenkov and Lamberti~\cite{buladir, bula} aiming at estimates via the Lebesgue measure of the symmetric difference of the open sets.

Namely, we consider families of open sets which are locally  subgraphs of functions of class $C^{m-1,1}$. We require that the `atlas'
${\mathcal{A}}$, with the help of which such boundaries are described, is fixed and we consider the class
$C^{m-1,1}_M({\mathcal{A}})$ of open sets for which the behavior of the derivatives of the functions describing the boundaries is controlled by a fixed constant $M>0$ (see Definition~\ref{class}).

Let $\varphi_n[\Omega ]$, $n\in {\mathbb{N}}$,  denote an orthonormal sequence of eigenfunctions corresponding
to the eigenvalues $\lambda _n[\Omega ]$.
In Corollary~\ref{unifbounded} we prove that if ${\mathfrak{A}}$ is a family of open sets of class $C^{m-1,1}_M({\mathcal{A}})$ such that for some $2<p\le \infty $
\begin{equation}
\label{unifbound} \sup_{\Omega \in
{\mathfrak{A}}}\|\varphi_n[\Omega] \|_{W^{m,p}(\Omega )} <\infty ,
\end{equation}
for all $n\in {\mathbb{N}}$,
then  for each $n\in\mathbb{N}$ there exists $c_n>0$ such
that for both Dirichlet and Neumann boundary conditions
\begin{equation}
\label{intro1}
|\lambda_n[\Omega_1]-\lambda_n[\Omega_2] |\le c_n |\Omega_1
\vartriangle \Omega_2 |^{1-\frac{2}{p}},
\end{equation}
for all $\Omega_1,\Omega_2\in {\mathfrak{A}}$ satisfying $ |\Omega_1
\vartriangle \Omega_2 |< c_n^{-1}$, where   $|\Omega_1
\vartriangle \Omega_2 |$ is the Lebesgue measure of the symmetric difference $ \Omega_1
\vartriangle \Omega_2 $.

If $\Omega_1$  is fixed and $\Omega _2\subset \Omega_1$ then  in the case of Dirichlet boundary conditions the assumptions of Corollary~\ref{unifbounded}
can be weakened. In fact, in Corollary~\ref{sharpinc} we prove that if $\Omega_1$ is of class $C^{m-1,1}_M({\mathcal{A}} ) $ and, for some $2<p\le \infty $,  $\varphi_n[\Omega_1]\in W^{m,p}(\Omega_1 )$ for all $n\in {\mathbb{N}}$, then
 for each $n\in\mathbb{N}$ there exists $c_n>0$ such
that for Dirichlet boundary conditions
 \begin{equation}
 \label{intro2}
 \lambda_{n}[\Omega_1]\le \lambda_{n}[\Omega_2]\le \lambda_{n}[\Omega_1]+c_n|\Omega_1\setminus \Omega_2|^{1-\frac{2}{p}},
 \end{equation}
for all $\Omega _2$ of class $C^{m-1,1}_M({\mathcal{A}} ) $ satisfying $\Omega_2\subset \Omega_1$ and
$|\Omega_1\setminus \Omega_2|\le c_n^{-1}$. (In this case there are no assumptions on the eigenfunctions $\varphi_n[\Omega_2]$.)

In Section 5 we also prove that, in general, the exponent $1-2/p$ in (\ref{intro1}) and (\ref{intro2}) cannot be replaced by a larger one.

If the coefficients $A_{\alpha \beta }$ are of class $C^m$ and the open sets are of class $C^{2m}$, condition (\ref{unifbound}) is satisfied with $p=\infty $. It follows that
for each $n\in\mathbb{N}$ there exists $c_n>0$ such that
\begin{equation}
\label{intro3}
|\lambda_n[\Omega_1]-\lambda_n[\Omega_2] |\le c_n |\Omega_1
\vartriangle \Omega_2 |,
\end{equation}
for all $\Omega_1,\Omega_2$ of class $C^{2m}_M({\mathcal {A}})$ satisfying $|\Omega_1\vartriangle \Omega_2|<c_n^{-1}$. See Corollary~\ref{smoothest}.

The case $m=1$ was considered in \cite{buladir, bula}. As in \cite{buladir, bula},
the proof of  our estimates is based on  the  general spectral stability theorem \cite[Thm.~3.2]{bula}.
In order to apply that theorem we  construct   linear operators  ${\mathcal{T}}_{{\mathcal{D}}}: W^{m,2}_0(\Omega _1)\to W^{m,2}_0(\Omega _2)$,
${\mathcal{T}}_{{\mathcal{N}}}: W^{m,2}(\Omega _1)\to W^{m,2}(\Omega _2)$
 possessing a number of  special properties. These operators serve as `transition operators'
 for Dirichlet and Neumann  boundary conditions respectively, as required by the general spectral stability theorem.
   We point out that the construction of such transition operators  for $m >1$ is rather sofisticated  and a straightforward extension to the case $m>1$ of the techniques used in \cite{buladir} for $m=1$ is not possible (see the beginning of Section \ref{section3} for details).

We note that in \cite{bulahigh} we proved spectral stability estimates expressed in terms of  so-called `atlas' distance introduced in
 \cite[Definition 5.1]{bulahigh}
 and of the Hausdorff  distance of the boundaries of $ \Omega_1$ and $\Omega_2$. In that case we  considered classes of open sets with boundaries admitting  arbitrarily strong degenerations and we did not require  any  summability  assumption on the eigenfunctions and their gradients. However, as we  pointed out in  \cite[Example 8.1]{buladir}, using the Lebesgue measure of $\Omega_1\vartriangle \Omega_2$ as we do here, allows to obtain better
 estimates.

\section{Preliminaries and notation}
\label{section2}

Let $N, m \in \mathbb{N}$ and $\Omega$ be an open set in $\mathbb{R}^N$. Let ${\mathbb{N}}_0^N$ be the set of all multi-indices $\alpha =(\alpha_1, \dots , \alpha_N)$ and $|\alpha |=\alpha_1+\dots +\alpha_N$ be their lengths. Here ${\mathbb{N}}_0={\mathbb{N}}\cup\{0\} $.
By $W^{m,p}(\Omega )$, $1\le p \le \infty $,  we denote the Sobolev space  of
all  complex-valued functions $u$ in $L^p(\Omega )$, which
have all weak derivatives $D^{\alpha }u$ up to order $m$ in $L^p(\Omega ) $,
endowed with the norm
\begin{equation}
\label{norm}
\|u \|_{W^{m,p}(\Omega )}=\sum_{|\alpha|\le m}\|D^{\alpha }u \|_{L^{p}(\Omega )}.
\end{equation}
If $1\le p<\infty $, then by  $W^{m,p}_0(\Omega )$ we denote
the closure in $W^{m,p}(\Omega )$ of the space of all $C^{\infty}$-functions
with compact support in $\Omega$. For open sets $\Omega$ under
consideration a function belongs to $W^{m,p}_0(\Omega )$ if and
only if its extension by zero outside $\Omega$ belongs to
$W^{m,p}({\mathbb{R}}^N)$. By
$W^{m,\infty}_0(\Omega)$ we denote the space of all functions in
$W^{m,\infty}(\Omega )$ whose extension by zero outside $\Omega$ belongs
to
$W^{m,\infty}({\mathbb{R}}^N)$, which is wider than the closure in
$W^{m,\infty}(\Omega )$ of the space of all $C^{\infty}$-functions
with compact support in $\Omega$.

Let $\hat m$ be the number of the
multi-indices $\alpha\in {\mathbb{N}}_0^N$
with $|\alpha | = m$.
For all $\alpha ,\beta \in
{\mathbb{N}}_0^N$
such that $|\alpha |=|\beta |=m$, let
$A_{\alpha \beta }$ be bounded measurable real-valued
functions defined on $\Omega$ such that
$A_{\alpha \beta }=
A_{\beta \alpha }$ and  for some $\theta >0$
\begin{equation}
\label{elp}
\theta^{-1} |\xi |^2\le   \sum_{|\alpha |=|\beta |=m}A_{\alpha \beta }(x)\xi_{\alpha }\xi_{\beta }  \le \theta |\xi |^2
\end{equation}
for all $x\in \Omega $, $\xi =(\xi_{\alpha })_{|\alpha |=m}\in\mathbb{R}^{\hat m}$.

Let $V(\Omega )$ be a closed subspace of  $W^{m,2}(\Omega )$ containing  $W^{m,2}_0(\Omega )$. We consider the
following eigenvalue problem

\begin{equation}
\label{mainprobl}
\int_{\Omega}\sum_{|\alpha |=|\beta |=m}
A_{\alpha \beta }D^{\alpha }u D^{\beta }\bar v dx=\lambda \int_{\Omega}u \bar v dx,
\end{equation}
for all test functions $v\in V (\Omega )$, in the unknowns
$u\in V (\Omega )$ (the eigenfunctions) and $\lambda \in\mathbb{R}$
(the eigenvalues).

As is well-known, problem (\ref{mainprobl}) is the weak formulation of the eigenvalue problem for the operator
$H$ in (\ref{classic}) subject to suitable homogeneous boundary conditions: the choice of $V(\Omega )$ corresponds to the choice of the boundary conditions (see {\it e.g.,} Ne\v{c}as \cite{Ne67}).

We set
\begin{equation}
Q_{\Omega }(u,v)=
\int_{\Omega }
\sum_{|\alpha |=|\beta |=m} A_{\alpha \beta }D^{\alpha }u
D^{\beta }\bar vdx,\ \ \ \ \ \ \ Q_{\Omega }(u)=Q_{\Omega }(u,u)  ,
\end{equation}
 for all $u,v\in W^{m,2}(\Omega ) $.

We assume that the  embedding $V(\Omega )\subset W^{m-1,2}(\Omega )$ is compact. Then one can prove that the restriction to $V(\Omega )$
of the quadratic form $Q_{\Omega }$ is closed,  hence the eigenvalues of equation (\ref{mainprobl}) coincide with the eigenvalues of
a suitable operator $H_{V(\Omega ) }$ canonically associated with $Q_{\Omega }$ and $V(\Omega)$. Since, in particular, the embedding $V(\Omega )\subset L^2(\Omega )$ is compact, $H_{V(\Omega ) }$ has compact resolvent and  the following theorem holds (see \cite[Thm.~2.1]{bulahigh} for a detailed proof).

\begin{thm}
\label{setting}
Let $\Omega$ be an open set in $\mathbb{R}^N$.
Let $m\in\mathbb{N}$, $\theta >0$ and, for all $\alpha ,\beta \in
 {\mathbb{N}}_0^{N} $
such that $|\alpha |=|\beta |=m$, let
$A_{\alpha \beta }$ be bounded measurable real-valued
functions defined on $\Omega$,  satisfying
$A_{\alpha \beta }= A_{ \beta\alpha  }$ and condition (\ref{elp}).

Let $V (\Omega )$ be a closed subspace of  $W^{m,2}(\Omega )$ containing  $W^{m,2}_0(\Omega )$  and such that the embedding
$V(\Omega )\subset W^{m-1,2}(\Omega )$ is compact.

Then there exists a non-negative self-adjoint linear operator
$H_{V(\Omega )}$ on $L^2(\Omega )$ with compact resolvent,
 such that ${\rm Dom}(H^{1/2}_{V(\Omega )})=
V (\Omega )$ and
\begin{equation}
\label{setting0}
<H^{1/2}_{V(\Omega )}u,H_{V(\Omega )}^{1/2} v >_{L^2(\Omega )}=
Q_{\Omega }(u,v)
,
\end{equation}
for all  $u,v \in V(\Omega )$.
Moreover, the eigenvalues of equation
(\ref{mainprobl})
coincide with the eigenvalues
$\lambda_n [H_{V(\Omega )} ]$ of
$H_{V(\Omega )} $ and
\begin{equation}
\label{setting1}
\lambda_n [H_{V(\Omega )}]
=\inf_{\substack{{\mathcal{L}}\subset   V (\Omega ) \\ {\rm dim}\, {\mathcal{L}}=n}}\sup_{\substack{u\in {\mathcal{L}} \\ u\ne 0}}
\frac{Q_{\Omega }(u)}{\|u \|_{L^2(\Omega )}^2}
 ,
\end{equation}
where the infimum  is taken with respect to all subspaces ${\mathcal{L}}$ of $V(\Omega )$ of dimension $n$.
\end{thm}

Note that the compactness of the embedding  $V(\Omega )\subset W^{m-1,2}(\Omega )$ can be deduced by the compactness of the embedding $V(\Omega )\subset L^2(\Omega )$ under some further assumptions on $\Omega$. Assume that $\Omega$ is such that
 for any $\epsilon >0$ there exists $c(\epsilon )>0$ such that the following inequality holds:
 $$
 \| u \|_{W^{m-1,2}(\Omega )}\le c(\epsilon )\| u\|_{L^2(\Omega )}+\epsilon \sum_{|\alpha |=m}\| D^{\alpha }u \|_{L^2(\Omega )}.
 $$
(This inequality holds in particular if $\Omega $ has a quasi-continuous boundary, see Burenkov~\cite[Thm.~6, p.~160]{bu}.) Then
the compactness of the embedding
 $V(\Omega )\subset L^2(\Omega )$ is equivalent to the compactness of the embedding $V(\Omega )\subset W^{m-1,2}(\Omega )$, see Burenkov~\cite[Lemma~13, p.~172]{bu} for details.

In this paper we are interested in the cases $V(\Omega )=W^{m,2}_0(\Omega )$ and $V(\Omega )=W^{m,2}(\Omega )$ which correspond to Dirichlet and Neumann boundary conditions respectively.

\begin{defn}\label{defset} Let $\Omega$ be an open set in ${\mathbb{R}}^N$.
Let $m\in\mathbb{N}$, $\theta >0$ and, for all $\alpha ,\beta \in
 {\mathbb{N}}_0^{N} $
such that $|\alpha |=|\beta |=m$, let
$A_{\alpha \beta }$ be bounded measurable real-valued
functions defined on $\Omega$,  satisfying
$A_{\alpha \beta }= A_{ \beta\alpha  }$ and condition (\ref{elp}).

If the embedding $W^{m,2}_0(\Omega )\subset  W^{m-1,2}(\Omega )$ is compact, we set
$$
\lambda_{n,{\mathcal {D}}}[\Omega ]=\lambda_n[H_{W^{m,2}_0(\Omega )}  ].
$$

If the embedding $W^{m,2}(\Omega )\subset  W^{m-1,2}(\Omega )$ is compact,
we set
$$
\lambda_{n,{\mathcal {N}}}[\Omega ]=\lambda_n[H_{W^{m,2}(\Omega )}  ].
$$

The numbers $\lambda_{n,{\mathcal {D}}}[\Omega ]$, $\lambda_{n,{\mathcal {N}}}[\Omega ]$   are called the Dirichlet eigenvalues, Neumann eigenvalues
respectively,  of  operator (\ref{classic}).
\end{defn}

\begin{rem}
\label{sufcomp}
If $\Omega $ is such that the embedding $W^{1,2}_0(\Omega )\subset L^2(\Omega )$ is compact (for instance, if $\Omega $ is an arbitrary open set with  finite Lebesgue measure), then also the embedding $ W^{m,2}_0(\Omega ) $ $  \subset W^{m-1,2}(\Omega )$ is compact and the Dirichlet eigenvalues are well-defined.

If $\Omega $ is such that the embedding $W^{1,2}(\Omega )\subset L^2(\Omega )$ is compact (for instance, if $\Omega $
has a continuous boundary, see Definition~\ref{class}), then the embedding $W^{m,2}(\Omega )\subset W^{m-1,2}(\Omega )$ is compact and
the Neumann eigenvalues are well-defined.
\end{rem}

In the next sections we shall study the variation of $\lambda_{n,{\mathcal{D}}}[\Omega ]$ and $\lambda_{n,{\mathcal{N}}}[\Omega ]$ upon variation of $\Omega$ in suitable classes of open sets defined below.

For any set $V$ in ${\mathbb{R}}^N$ and $\delta >0$ we denote by $V_{\delta }$ the set $\{x\in V:\ d(x, \partial \Omega )>\delta \}$. Moreover,   as in \cite{buda}, by a cuboid we mean  any rotation of a rectangular parallelepiped in ${\mathbb{R}}^N$.\\

\begin{defn}
\label{class}

Let $ \rho >0$, $s,s'\in\mathbb{N}$, $s'\le s$
and  $\{V_j\}_{j=1}^s$ be a family of bounded open cuboids  and
$\{r_j\}_{j=1}^{s} $ be a family of rotations in ${\mathbb{R}}^N $.

We say that ${\mathcal{A}}= (  \rho , s,s', \{V_j\}_{j=1}^s, \{r_j\}_{j=1}^{s} ) $ is an atlas in ${\mathbb{R}}^N$ with the parameters
$\rho , s,s', \{V_j\}_{j=1}^s, \{r_j\}_{j=1}^{s}$, briefly an atlas in ${\mathbb{R}}^N$.

We denote by $C( {\mathcal{A}}   )$ the family of all open sets $\Omega $ in ${\mathbb{R}}^N$
satisfying the following properties:

(i) $ \Omega\subset \bigcup\limits_{j=1}^s(V_j)_{\rho}$ and $(V_j)_\rho\cap\Omega\ne\emptyset;$

(ii) $V_j\cap\partial \Omega\ne\emptyset$ for $j=1,\dots s'$, $ V_j\cap \partial\Omega =\emptyset$ for $s'<j\le s$;

(iii) for $j=1,...,s$
$$
r_j(V_j)=\{\,x\in \mathbb{R}^N:~a_{ij}<x_i<b_{ij}, \,i=1,....,N\}
$$

\noindent and

$$
r_j(\Omega\cap V_j)=\{x\in\mathbb{R}^N:~a_{Nj}<x_N<g_{j}(\bar
x),~\bar x\in W_j\},$$

\noindent where $\bar x=(x_1,...,x_{N-1})$, $W_j=\{\bar
x\in\mathbb{R}^{N-1}:~a_{ij}<x_i<b_{ij},\,i=1,...,N-1\}$
and $g_j$ is a continuous function defined on $\overline {W}_j$ (it is meant that if $s'<j\le s$ then $g_j(\bar x)=b_{Nj}$ for all $\bar x\in \overline{W}_j$); moreover for $j=1,\dots ,s'$
$$
a_{Nj}+\rho\le g_j(\bar x)\le b_{Nj}-\rho ,$$

\noindent for all $\bar x\in \overline{W}_j$.

We say that an open set $\Omega$ in ${\mathbb{R}}^N$ is an open set with a continuous boundary if $\Omega $ is of class
$C( {\mathcal{A}}   ) $ for some atlas $ {\mathcal{A}}   $.

Let $m\in {\mathbb{N}}, M>0$. We say that an open set $\Omega$ is of class $C^m_M( {\mathcal{A}}   )$, $C^{m-1,1}_M( {\mathcal{A}}   )  $  if $\Omega $ is of class $C( {\mathcal{A}}   )$
and all the functions $g_j$ in (iii) are of class $C^m(\overline{W}_j)$, $C^{m-1,1}(\overline{W}_j)$ with
$$
| g_j |_{c^m(\overline {W}_j)} =  \sum_{1\le |\alpha |\le m}\| D^{\alpha }g_j\|_{L^{\infty }(\overline{W}_j)}\le M,
$$
$$
\ |g_j|_{  c^{m-1,1}(\overline {W}_j)   }=|g_j|_{c^{m-1}(\overline {W}_j)}+\sum_{|\alpha |=m-1} \sup_{\substack{\bar x,\bar y\in {\overline {W}}_j\\ \bar x\ne \bar y}}\frac{|D^{\alpha }g_j(\bar x)-D^{\alpha }g_j(\bar y)|}{|\bar x-\bar y|}\le M
$$
respectively\footnote{Note that as customary $\| g_j\|_{C^m  (\overline{W}_j)  }=\| g_j\|_{L^{\infty }(W_j)} + | g_j |_{c^m(\overline {W}_j)}$  and
$\| g_j\|_{C^{m-1,1}  (\overline{W}_j)  }=\| g_j\|_{L^{\infty }(W_j)} + | g_j |_{c^{m-1,1}(\overline {W}_j)}$.}.

We say that an open set $\Omega$ in ${\mathbb{R}}^N$ is an open set of class $C^m$, $C^{m-1,1}$  if $\Omega $ is of class
$C^m_M( {\mathcal{A}}   ) $, $C^{m-1,1}_M( {\mathcal{A}}   ) $ respectively, for some atlas $ {\mathcal{A}}   $ and some $M>0$.

\end{defn}

\section{A pre-transition operator for higher order So\-bo\-lev spaces}
\label{section3}

The aim of this section is proving the following theorem.

\begin{thm}\label{pretran}
Let
${\mathcal{A}}= (  \rho , s,s', \{V_j\}_{j=1}^s, \{r_j\}_{j=1}^{s} ) $ be an atlas in ${\mathbb{R}}^N$, $m\in {\mathbb{N}},$ $ M>0$. Let $\Omega_1,\Omega_2\in C^{m-1,1}_M( {\mathcal{A}}   )  $. For all  $m\in {\mathbb{N}}$, $1\le p\le \infty $
there exist linear maps
$$ {\mathcal{T}}_{{\mathcal{D}}}:W^{m,p}_0(\Omega_1)\to W^{m,p}_0(\Omega_2)
\ \ {and}\ \ {\mathcal{T}}_{{\mathcal{N}}}:W^{m,p}(\Omega_1)\to W^{m,p}(\Omega_2),
$$
with the following properties:
\begin{enumerate}
\item[(i)]
there exists $C_1>0$ depending only on $ {\mathcal{A}}, m, M, p$ such that
$
\| {\mathcal{T}}_{{\mathcal{D}}}\|, \| {\mathcal{T}}_{{\mathcal{N}}}\| \le C_1.
$
\item[(ii)] there exists $C_2>0$ depending only on $ {\mathcal{A}}$, and an open set $\Omega_3\subset \Omega_1\cap\Omega_2$ such that
\begin{equation}
\label{measure}
|\Omega _1\setminus \Omega_3|,\ |\Omega _2\setminus \Omega_3|\le C_2| \Omega_1 \vartriangle \Omega_2  |,
\end{equation}
and such that
\begin{equation}
\label{identity}
{\mathcal{T}}_{{\mathcal{D}}}[u](x)=u(x),\ \ \ \ {\mathcal{T}}_{{\mathcal{N}}}[v](x)=v(x),
\end{equation}
for all $u\in  W^{m,p}_0(\Omega_1), v\in W^{m,p}(\Omega_2)$,  $x\in \Omega_3$.
\end{enumerate}
\end{thm}

For $m=1$ Theorem~{\ref{pretran}} was proved in \cite{buladir, bula}. We note that the construction of the operator ${\mathcal{T}}_{{\mathcal{D}}}$ in \cite{buladir} does not have a straightforward generalization to the case $m>1$. A more or less straightforward generalization of a crucial step in the construction in \cite{buladir} is given in Lemma~\ref{graf}, where a special transformation $\Phi_c$ is defined depending on a positive constant $c$. Importantly, the derivatives of $\Phi_c$ of order greater than one have singularities. This leads to singularities when applying the chain rule to compositions $v(\Phi_c)$, which are explicitly written out in the first summand of the right-hand side of formula (\ref{monfor}). In order to overcome this difficulty, in Lemma~\ref{paste} we construct a linear map $T$, given by formula (\ref{paste2}), with appropriately chosen parameters $\delta_k$, $c_k$ which allow to control the effect of singularities and ensure the boundedness of $T$. The proof of the boundedness of $T$ is based on the one-dimensional Lemma~\ref{ext}. Finally, in the proof of Theorem~\ref{pretran} local transformations of such type are pasted together.

\begin{lem}
\label{graf} Let  $W$ be a bounded convex open set in
$\mathbb{R}^{N-1}$. Let $m\in \mathbb{N}$, $a\in \mathbb{R}$, $D_1>D_2>a$ and $g_1,g_2\in C^{m-1,1}(\overline{W})$ be such that
\begin{equation}
\label{graf1}
D_2<g_2(\bar x),\  g_1(\bar x)<D_1 ,
\end{equation}
for all $\bar x\in \overline{W}$. Let $\delta =
\frac{D_2-a}{2(D_1-D_2)}$ and $c\geq \frac{1}{\delta } $. Let
\begin{eqnarray*}
& & g_3(\bar x)  = g_2(\bar x)-\delta (g_1(\bar x)-g_2(\bar x))^+\\
& & g_{1,c}(\bar x)  = g_2(\bar x)+c\delta (g_1(\bar x)-g_2(\bar x))^+,
\end{eqnarray*}
for all $\bar x\in \overline{W}$, and let
\begin{eqnarray}
\label{graf2}
& &{\mathcal{O}}_k= \left\{(\bar x,x_N):\ \bar x\in W,\ a<x_N<g_k(\bar x) \right\},\ \ \ \ \ k=1,2,3,
\nonumber\\
& &{\mathcal{O}}_{1,c}= \left\{(\bar x,x_N):\ \bar x\in W,\ a<x_N<g_{1,c}(\bar x) \right\}.
\end{eqnarray}
Let $\Phi_{c}$ be the map of $\overline{{\mathcal{O}}}_2$ into $\overline{{\mathcal{O}}}_{1c}$
defined by
\begin{equation}
\label{graf3}
\Phi_{c}(x)= (\bar x, x_N+ch(x)),\ \ x\in {\overline {\mathcal {O}}_2}
\end{equation}
where
\begin{equation}\label{acca}
h(x)= \left\{
\begin{array}{ll}
0, & {\rm if}\ x\in  {\overline{{\mathcal{O}}}}_{3 },\\
\frac{(x_N-g_3(\bar x))^{m+1}}{\delta^m(g_1(\bar x)-g_2(\bar x))^m},
 & {\rm if}\ x\in {\overline{{\mathcal{O}}}}_2\setminus   {\overline{{\mathcal{O}}}}_{3},
\end{array}
\right.
\end{equation}
Then the following statements hold:
\begin{itemize}
\item[(i)] $\emptyset \ne {\mathcal{O}}_{3}\subset {\mathcal{O}}_{2}$; $  {\mathcal{O}}_{1}, {\mathcal{O}}_{2}, {\mathcal{O}}_{3}\subset {\mathcal{O}}_{1,c}$\ \ and\ \   $
|{\mathcal{O}}_{1,c}\setminus {\mathcal{O}}_{2}|= c\delta  |{\mathcal{O}}_{1}\setminus {\mathcal{O}}_{2}|=c |{\mathcal{O}}_{2}\setminus {\mathcal{O}}_{3}| $;
\item[(ii)] $\Phi_{c}$ is a bijection of
$\overline{{\mathcal{O}}}_2$ onto $\overline{{\mathcal{O}}}_{1,c}$,\  $\Phi_{c}(\partial {\mathcal{O}}_{2})=\partial
{\mathcal{O}}_{1,c}$, $\Phi_{c}\in C^{m-1,1}_{loc}({\mathcal{O}}_{2} )\cap {\rm Lip ({\mathcal{O}}_{2} )}$, and $\Phi_c(x)=x$ for all $x\in {\mathcal{O}}_3$;
\item[(iii)]
there exists $M>0$ depending only on $N, m, a, D_1, D_2, \|g_1\|_{c^{m-1,1}(\overline{ W})}$ and $\|g_2\|_{c^{m-1,1}(\overline{ W})}$ such that
for all $\alpha \in {\mathbb{N}}_0^N$ with $|\alpha |\le m$
$$
\|h^{|\alpha |-1}D^{\alpha }h\|_{ L^{\infty }( {{\mathcal{O}}}_2\setminus   {\overline{{\mathcal{O}}}}_{3})    }\le
M.
$$
\end{itemize}
\end{lem}

{\bf Proof.} We  note that if $\bar x\in \overline{W}$ and $g_2(\bar x)\geq g_1(\bar x)$ then $g_3(\bar x)= g_{1,c}(\bar x)=g_2(\bar x)$; viceversa, if $g_2(\bar x)< g_1(\bar x)$, since $c\delta \geq 1$ it follows that $g_3(\bar x)< g_2(\bar x)< g_1(\bar x) <g_{1,c}(\bar x)$.
 In particular, $ {\mathcal{O}}_{3}\subset {\mathcal{O}}_{2}$ and  $  {\mathcal{O}}_{1}, {\mathcal{O}}_{2}, {\mathcal{O}}_{3}\subset {\mathcal{O}}_{1,c}$; moreover, if $(\bar x, x_N)\in {\overline{{\mathcal{O}}}}_2\setminus   {\overline{{\mathcal{O}}}}_{3}$ then $g_2(\bar x)<g_1(\bar x)$, hence $\Phi_c$ is well defined. Since $\delta < (D_2-a)/(D_1-D_2)$ then
 $a<g_3(\bar x) $ for all $\bar x\in \overline {W} $, hence ${\mathcal{O}}_{3}\ne \emptyset$. Moreover, we note that
 $$
 g_{1,c}(\bar x)-g_2(\bar x)=c\delta (g_1(\bar x)-g_2(\bar x))^+=c(g_2(\bar x)-g_3(\bar x)),
 $$
hence the equalities in statement {\it (i)} follow.

Statement {\it (iii)}  follows by standard calculus.

 We now prove statement {\it (ii)}. By using the same argument as in \cite[Lemma~4.1]{buladir}
 one can  prove that $\Phi_c\in {\rm Lip} ({\mathcal{O}}_{2} )$. Moreover, it is obvious that $\Phi_c$ is a bijection
 of $\overline {\mathcal{O}}_{2}  $ onto $\overline {\mathcal{O}}_{1,c}  $ and $\Phi_{c}(\partial {\mathcal{O}}_{2})=\partial
{\mathcal{O}}_{1,c}$.

It remains to prove that $\Phi_c\in C^{m-1,1}_{loc}({\mathcal{O}}_2)$. Clearly $\Phi_c$ is of class $C^{m-1,1}_{loc}$ on the open sets  ${\mathcal{O}}_{3}$ and
 $ {\mathcal{O}}_{2}\setminus {\overline{\mathcal{O}}}_{3}$. We now prove that $\Phi_c$ is of class $C^{m-1}$ in a neighborhood of any point of
 $ {\mathcal{O}}_{2}\cap \partial {\mathcal{O}}_{3} $. It clearly suffices to do so for $(\Phi_c)_N$. Let $(\bar y, y_N)\in {\mathcal{O}}_{2}\cap \partial {\mathcal{O}}_{3} $. Then
$y_N=g_3(\bar y)<g_2(\bar y)< g_1(\bar y)$ and by continuity there exists an open  neighborhood $U$ of $(\bar y, y_N)$ contained in ${\mathcal{O}}_{2}$ such that
$g_3(\bar x)<g_2(\bar x)< g_1(\bar x)$ for all $(\bar x, x_N)\in U$.

Consider  the functions $\varphi_1(\bar x, x_N)= x_N$ and $\varphi_2(\bar x ,x_N)= x_N+\frac{c|x_N-g_3(\bar x)|^{m+1}}{\delta^m(g_1(\bar x)-g_2(\bar x))^m}$ for all $(\bar x , x_N)\in U$. Clearly $\varphi_1, \varphi_2\in C^{m-1,1}(U) $ and $D^{\alpha }\varphi _1
= D^{\alpha }\varphi_2 $ on $U \cap \partial {\mathcal{O}}_{3}$  for all $\alpha \in {\mathbb{N}}_0^N$ with $|\alpha |\le m-1$. Since
$(\Phi _c)_N= \varphi_1$ on $U \cap \overline {\mathcal{O}}_{3}$ and $(\Phi _c)_N= \varphi_2$ on $U\setminus \overline {\mathcal{O}}_{3} $ it follows that
$(\Phi_c)_N \in C^{m-1}(U)$.  Moreover, since $D^{\alpha }\varphi _1
= D^{\alpha }\varphi_2 $  on the graph of $g_3$ for all $|\alpha |=m-1$, it follows that $D^{\alpha }(\Phi _c)_N $, hence $D^{\alpha  }\Phi _c $   is locally  Lipschitz continuous on ${\mathcal{O}}_2$  for all $|\alpha |=m-1$.\hfill $\Box$\\

\begin{lem}\label{mon-1}
Let the assumptions of Lemma~\ref{graf} hold. If $v\in W^{m,1}_{loc}({\mathcal{O}}_{1,c} ) $ then $v\circ \Phi_c\in W^{m,1}_{loc}({\mathcal{O}}_{2} )$ and for each $\alpha\in {\mathbb{N}}_0^N$ with $1\le |\alpha |\le m$
\begin{eqnarray}
\label{monfor}\lefteqn{
D^{\alpha }(v( \Phi_c ))(x)} \\ & & =
\sum_{1\le |\beta |< |\alpha |}\frac{(D^{\beta }v)(\Phi_c(x))}{h(x)^{|\alpha |-|\beta |} }\sum_{r=1}^{| \beta  |} b_{\beta ,r }(x)c^{r}
+\sum_{|\beta |= |\alpha |}(D^{\beta }v)(\Phi_c(x))\sum_{r=0}^{| \beta  |} b_{\beta ,r }(x)c^{r}
,\nonumber
\end{eqnarray}
for all $x\in {\mathcal{O}}_{2} \setminus \overline {{\mathcal{O}}}_3$, where   $b_{\beta ,r  }$ are bounded continuous functions independent of $c$. Moreover, there exists $M>0$   depending only
on $N, m, a,$ $ D_1, D_2,  \|g_1\|_{c^{m-1,1}(\overline{ W})}$ and $\|g_2\|_{c^{m-1,1}(\overline{ W})}$ such that all functions  $b_{\beta ,r} $ in (\ref{monfor}) satisfy the inequality
$$
\| b_{\beta , r}\|_{L^{\infty }({\mathcal{O}}_2 \setminus \overline {{\mathcal{O}}}_3)} \le M.
$$

\end{lem}

{\bf Proof.}
 If $\phi =(\phi_1, \dots , \phi_N)$ is a map of ${\mathcal{O}}_2$ to ${\mathcal{O}}_{1,c}$ of class $C^{m-1,1}_{loc}$ then   $v\circ \phi \in W^{m,1}_{loc}({\mathcal{O}}_{2} )$ for all $v\in W^{m,1}_{loc}({\mathcal{O}}_{1,c} ) $. Moreover, by the chain rule $D^\alpha(v(\phi))$ is a linear combination of the functions
 \begin{equation}\label{mon0} (D^\beta v)(\phi)D^{\nu_{i_1}}\phi_{i_1}\cdots D^{\nu_{i_k}}\phi_{i_k}
  \end{equation}
with natural coefficients depending only on $ \alpha, \beta , \nu_{i_1},\dots , \nu_{i_k}$, where  $1\le|\beta|\le|\alpha|$,
  $k= |\beta |$, $i_{1},\dots i_{k}\in \{1, \dots , N\}$,  $\nu_{i_1},\dots,\nu_{i_k}\in{\mathbb{N}}_0^N$,  and
\begin{equation}
\label{coefficients} |\nu_{i_1}|+\cdots+|\nu_{i_k}|=|\alpha |,~~|\nu_{i_1}|,\dots,|\nu_{i_k}|\ge 1\,.
\end{equation}

In particular if $\phi =\Phi_c$
 then $\phi_i(x)=x_i$ for all $i=1,\dots , N-1$, and
  $\phi_N(x)=x_N+ch(x) $ for all $x\in {\mathcal{O}}_{2} $, where $h$ defined by (\ref{acca}).
If $i_{1}, \dots , i_{k}\in \{1, \dots , N-1 \}$ then among the  functions in ({\ref{mon0}})  we can consider only those with
$
\nu_{i_1}=e_{i_1},\dots,\nu_{i_k}= {e_{i_k}},
$ (here $e_1, \dots , e_N$ denotes the canonical basis in ${\mathbb{R}}^N$)
in which case $|\beta |=|\alpha |$ by (\ref{coefficients}):\; thus, in this case we can consider only functions of the type
\begin{equation}
\label{monsim}
(D^\beta v)(\Phi _c)
\end{equation}
with   $|\beta |=|\alpha |$. The remaining functions correspond to the cases when at least one of the indices $i_s$ is $N$.
Assume that exactly $n$ of them are equal to $N$, then $\nu_{i_s}=e_{i_s}$ for the remaining $|\beta |-n$ of them. Thus, such
 functions are of the type
\begin{equation}\label{mon1} (D^\beta v)(\Phi _c)D^{\eta_{1}}(x_N+ch(x))\cdots D^{\eta_{n}} (x_N+ch(x))
\end{equation}
where  $1\le|\beta|\le|\alpha|$,
  $1\le n\le |\beta |$,  $\eta_{1},\dots,\eta_{n}\in{\mathbb{N}}_0^N$,  and
\begin{equation}
\label{coefficientseta} |\eta_{1}|+\cdots+|\eta_{n}|=|\alpha|-|\beta |+n,~~|\eta_{1}|,\dots,|\eta_{n}|\ge 1\,.
\end{equation}
The functions in (\ref{mon1}) are linear combinations of functions of the type
\begin{equation}\label{mon2} c^{\rho }(D^\beta v)(\Phi_c)D^{\xi_{1}}h\cdots D^{\xi_{\rho }} h
\end{equation}
with natural coefficients depending only on $\alpha , \beta , \rho , \xi_1, \dots , \xi_{\rho } $, where
\begin{equation}
\label{coefficientsrho}  1\le \rho \le n,\ \ |\xi_{1}|+\cdots+|\xi_{\rho}|=|\alpha|-|\beta |+\rho,~~|\xi_{1}|,\dots,|\xi_{n}|\ge 1\, ,
\end{equation}
and of functions of the type (\ref{monsim}) which correspond to the case  $\eta_1=\dots =\eta_n=e_N$, in which case
$| \alpha  |=| \beta | $ by (\ref{coefficientseta}).

By Lemma~\ref{graf} the functions  $b_{\xi_s}=h^{|\xi_s |-1}  D^{\xi_s}h$  are continuous, bounded and such that  $\| b_{\xi_s}\|_{L^{\infty }( {\mathcal{O}}_{2}\setminus {\overline{\mathcal{O}}}_{3} ) }\le M_1$ where $M_1$ is constant depending only on $N, m, a,$ $ D_1, D_2,  \|g_1\|_{c^{m-1,1}(\overline{ W})}$ and $\|g_2\|_{c^{m-1,1}(\overline{ W})}$. Hence
\begin{eqnarray}\label{crho}\lefteqn{c^{\rho }(D^\beta v)(\Phi_c)D^{\xi_{1}}h\cdots D^{\xi_{\rho }} h} \\  \nonumber & & = c^{\rho }(D^\beta v)(\Phi_c)b_{\xi_1 }\cdots b_{\xi _{\rho }}h^{\rho-(|\xi_1|+\dots +|\xi _{\rho }|)} = c^{\rho }(D^\beta v)(\Phi_c)b_{\xi_1 }\cdots b_{\xi _{\rho }}h^{|\beta |-|\alpha |}.\end{eqnarray}

Thus $D^\alpha(v(\Phi_c))$ is a linear combination with natural  coefficients depending only on  $\alpha , \beta , \rho , \xi_1, \dots , \xi_{\rho } $ of functions of the type (\ref{monsim}) with $|\beta |=|\alpha |$ and of the type  (\ref{crho}) with $1\le |\beta |\le |\alpha |$ and $1\le \rho\le |\beta |$ . Clearly, such functions can be arranged as in formula (\ref{monfor}).\hfill $\Box$

\begin{lem}\label{ext}
Let $\mu , s\in {\mathbb{N}}$, $s\le \mu$, $0<c_1<\dots <c_{\mu }$, $\gamma_1, \dots ,\gamma_{\mu }\in {\mathbb{R}}$ and
\begin{equation}
\label{van}
\sum_{k=1}^{\mu }\gamma_kc_k^{\sigma }=0,\ \ \ \sigma=1,\dots ,s.
\end{equation}
Moreover, let $-\infty <a<b<\infty $, $\eta \in C^1[a,b]$ and $\eta (x)> 0$, $\eta '(x)\geq 0$ for all $x\in (a,b)$.
Then
$$
\biggl\| \eta ^{-s}(x)\sum_{k=1}^{\mu }\gamma_kc_kf(x+c_k\eta (x))\biggr\|_{L^p(a,b)}\le C\| f^{(s)} \|_{L^p(a,b+c_{\mu }\eta (b))},
$$
for all $1\le p\le \infty$ and for all $f\in W^{s,p}(a,b+c_{\mu }\eta (b))$, where $$C=\sum_{k=1}^{\mu }|\gamma_k|c_k(c_k-c_1)^s.$$
\end{lem}
{\bf Proof.}
 If $f\in W^{s,p}(a,b+c_{\mu }\eta (b)))$ there exists a function $g$ equivalent to $f$ on $(a,b+c_{\mu }\eta (b))$ whose derivative $g^{(s-1)}$ is absolutely continuous on $[a,b+c_{\mu }\eta (b)]$.
By (\ref{van}) for all $x\in (a,b)$
\begin{eqnarray}& &
\sum_{k=1}^{\mu } \gamma_kc_kg(x+c_k\eta (x))=\sum_{k=1}^{\mu }\gamma_kc_k(g(x+c_k\eta (x))-g(x+c_1\eta (x)) )\nonumber \\ & &
=\eta (x)\sum_{k=1}^{\mu }\gamma_kc_k(c_k-c_1)\int_0^1g'(x+c_1\eta (x)+t_1(c_k-c_1)\eta (x))dt_1\nonumber \\ & &
=\eta (x)\sum_{k=1}^{\mu }\gamma_kc_k(c_k-c_1)\int_0^1g'(x+c_1\eta (x)+t_1(c_k-c_1)\eta (x))-g'(x+c_1\eta (x))dt_1\nonumber\\ & &
=\eta ^2(x)\sum_{k=1}^{\mu }\gamma_kc_k(c_k-c_1)^2\int_0^1\int_0^1t_1g''(x+c_1\eta (x)+t_1t_2(c_k-c_1)\eta (x))dt_2dt_1=\dots\nonumber \\ & &
 =\eta ^s(x)\sum_{k=1}^{\mu }\gamma_kc_k(c_k-c_1)^s\int_0^1\dots \int_0^1t_1\cdots t_{s-1}g^{(s)}(x+c_1\eta (x)+t_1\cdots t_s(c_k-c_1)\eta (x))dt_s\dots dt_1  .\nonumber \\
\end{eqnarray}
By Minkowski's inequality for integrals
\begin{eqnarray}\lefteqn{
\biggl\| \eta^{-s}(x)\sum_{k=1}^{\mu } \gamma_kc_kf(x+c_k\eta (x))  \biggr\|_{L^p(a,b)}=\biggl\| \eta ^{-s}(x)\sum_{k=1}^{\mu } \gamma_kc_kg(x+c_k\eta (x))  \biggr\|_{L^p(a,b)}} \nonumber \\  & & \le \sum_{k=1}^{\mu }|\gamma_k|c_k(c_k-c_1)^s\int_0^1\dots\int_0^1\| g^{(s)}(x+c_1\eta (x)+t_1\cdots t_s(c_k-c_1)\eta (x))\|_{L^p(a,b)}dt_s\dots dt_1 \, .\nonumber
\end{eqnarray}
Let $y=x+c_1\eta (x)+t_1\cdots t_s(c_k-c_1)\eta (x)$. Note that for all $x\in [a,b]$, $a\le y(x)\le b+c_{\mu }\eta (b)$ and $y'(x)\geq 1$. Hence
\begin{eqnarray}\lefteqn{
\| g^{(s)}(x+c_1\eta (x)+t_1\cdots t_s(c_k-c_1)\eta (x))\|_{L^p(a,b)}}\nonumber \\ & & \qquad\qquad\qquad\qquad \le \| g^{(s)}(y)\|_{L^p(a, b+c_{\mu }\eta (b))}=\| f^{(s)}\|_{L^p(a, b+c_{\mu }\eta (b))}
\end{eqnarray}
and the statement follows.\hfill $\Box$
\\

\begin{lem}
\label{paste} Let  $W, m, a, D_1,D_2, \delta , g_1,g_2$ be as in Lemma~\ref{graf}.
  Let $\delta_1, \dots , \delta_{m}\in {\mathbb{R}}$, $1/\delta \le c_1<\dots <c_{m}   $   be such that
\begin{equation}\label{paste1}\sum_{k=1}^{m}\delta_k=1,\ \ \ {\rm and}\  \ \
\sum_{k=1}^{m}\delta_kc_k^{\tau }=0,\ \ \ \ \tau =1, \dots , m-1.
\end{equation}
Let $g_3$, $g_{c_k}$, ${\mathcal{O}}_1$, ${\mathcal{O}}_2$, ${\mathcal{O}}_3$, ${\mathcal{O}}_{1,c_k}$, $\Phi_{c_k}$ be as in Lemma~\ref{graf} with
$c$ replaced by $c_k$ for all $k=1, \dots , m$.
Let $T$ be the linear map of $L^1_{loc }(W\times ]a,\infty [)$
to $L^1_{loc }( {\mathcal{O}}_2 )$ defined by
 \begin{equation}
 \label{paste2}
T[v](x)=\sum_{k=1}^{m}\delta_kv(\Phi_{c_k}(x)),
\end{equation}
for all $x\in {\mathcal{O}}_2$ and for all $v\in L^1_{loc }(W\times ]a,\infty [)$.

Then  for all $1\le p\le \infty $
\begin{equation}\label{paste3}T: W^{m,p}(W\times ]a,\infty [) \to
W^{m,p}({\mathcal{O}}_2) \end{equation}
and  there exists $C>0$ depending only on
  $N, m, p, a,     D_1, D_2,  \delta , c_m,  \|g_1\|_{c^{m-1,1}(\overline{ W})}$, $  \|g_2\|_{c^{m-1,1}(\overline{ W})}$ such that
  $\| T\| \le C$.

  Moreover, $T[v](x)=v(x) $ for all $v\in L^1_{loc }(W\times ]a,\infty [)$, $x\in {\mathcal{O}}_3$, and if $v= 0$ on ${\mathcal{O}}_1^{c} $ then
  $T[v]=0 $ on ${\mathcal{O}}_2^{c}$.
\end{lem}

{\bf Proof.}
First of all we recall that if $x\in {\mathcal{O}}_3$ then $\Phi_{c_k}(x)=x$ for all $k=1,\dots , m$. Thus by the first condition in (\ref{paste1})  it follows
that
\begin{equation}\label{paste4,5}
T[v](x)=\sum_{k=1}^{m}\delta_kv(\Phi_{c_k}(x))=\sum_{k=1}^{m}\delta_kv(x)=v(x),
\end{equation}
for all $v\in L^1_{loc }(W\times ]a,\infty [)$, $x\in {\mathcal{O}}_3$.

Let $v\in W^{m,p}(W\times ]a,\infty [)$. By Lemma~\ref{mon-1} $T[v]\in W^{m,1}_{loc}({\mathcal{O}}_2)$ and
\begin{equation}\label{paste5}
D^{\alpha }(  T[v](x))=\sum_{1\le |\beta |\le |\alpha |}\sum_{r=0}^{| \beta  |}\frac{ b_{\beta ,r }(x)}{h(x)^{|\alpha |-|\beta |}}\sum_{k=1}^{m}\delta _k c_k^{r} (D^{\beta }v)(\bar x, x_N +c_kh(x))  ,
\end{equation}
for all $|\alpha |=m$ and  for all $x\in {\mathcal{O}}_2\setminus \overline{{\mathcal{O}}}_3$, where $b_{\beta ,0}=0$ if $|\beta |<|\alpha |$.
We now estimate the $L^p$ norms of the summonds in the right-hand side of (\ref{paste5}).
We consider first the case $1\le p<\infty$, $|\beta |< |\alpha |$, $1\le r\le |\beta |$.  In this case we apply Lemma \ref{ext} with $f(x_N)=D^{\beta }v(\bar x , x_N ) $, $a=g_3(\bar x)$, $b= g_2(\bar x)   $,  $\eta (x_N)=h(\bar x, x_N)$,   $\mu =m, \gamma_k= \delta_kc_k^{r-1}, s=|\alpha |-|\beta|$.

Note that by (\ref{paste1})
$$
\sum_{k=1}^m\gamma_kc_k^{\sigma }  =
\sum_{k=1}^m\delta _kc_k^{\sigma +r-1 }=0,\ \ \ \sigma =1,\dots , |\alpha |-|\beta |.
$$
Indeed, $1\le r \le |\beta |$, hence $1\le \sigma +r-1\le |\alpha |-1\le m-1$ for all $\sigma =1,\dots , |\alpha |-|\beta |$.
Thus, condition (\ref{van}) is satisfied and by Lemma \ref{ext} we have
\begin{eqnarray}\label{cuno}\label{label6}\lefteqn{
\left(\int_W\int_{g_3(\bar x)}^{g_2(\bar x)}\left| \frac{1}{h(x)^{|\alpha |-|\beta |} } \sum_{k=1}^{m}  \delta_kc_k^{r} (D^{\beta }v)(\bar x, x_N +c_kh(x)) \right|^p dx_Nd\bar x\right)^{1/p}}\nonumber \\ & &
\le \tilde c_1\left(\int_W\int_{a}^{\infty }|D^{(\bar \beta ,\beta_N+|\alpha |-|\beta | )}v(\bar x, x_N)|^pdx_Nd\bar x\right)^{1/p}
\le \tilde c_1 \| v\|_{W^{m,p}(W\times ]a,\infty [)}\nonumber ,\\
\end{eqnarray}
where $\beta = (\bar \beta , \beta_N )$.
In the case $1\le p<\infty$, $|\beta |=|\alpha |$, $0\le r\le |\beta |$,  by a simple change of variables we obtain
\begin{eqnarray}\label{label7}\lefteqn{
\left(\int_W\int_{g_3(\bar x)}^{g_2(\bar x)}\left|  \sum_{k=0}^{m}\delta _kc_k^{r} \frac{(D^{\beta }v)(\bar x, x_N +c_kh(x))}{h(x)^{|\alpha |-|\beta |} } \right|^p dx_Nd\bar x\right)^{1/p}}\nonumber \\ & &
\le  \left(\int_W\int_{a}^{\infty  }\left|  \sum_{k=0}^{m} \delta _kc_k^{r}(D^{\beta  }v)(\bar x, x_N +c_kh(x)) \right|^p dx_Nd\bar x\right)^{1/p}
 \le \tilde c_2 \| v\|_{W^{m,p}(W\times ]a,\infty [)}.\nonumber \\
 \nonumber \\ & &
\end{eqnarray}

Thus by (\ref{paste5}), (\ref{label6}), (\ref{label7})
\begin{equation}\label{c3}
\|D^{\alpha } T[v]\|_{L^p({\mathcal{O}}_2\setminus {\mathcal{O}}_3)}\le \tilde c_3 \| v\|_{W^{m,p}(W\times ]a,\infty [)},
\end{equation}
for all $v\in W^{m,p}(W\times ]a,\infty [)$.
Clearly, by (\ref{paste4,5})
$$
\|D^{\alpha } T[v]\|_{L^p({\mathcal{O}}_3)}=\|D^{\alpha }v\|_{L^p({\mathcal{O}}_3)}\le \| v\|_{W^{m,p}(W\times ]a,\infty [)},
$$
Thus
\begin{equation}\label{c4}
\| T[v]\|_{W^{m,p}({\mathcal{O}}_2)}\le \tilde c_4 \| v\|_{W^{m,p}(W\times ]a,\infty [)}.
\end{equation}
for all $v\in W^{m,p}(W\times ]a,\infty [)$. In (\ref{cuno})-(\ref{c4}), $\tilde c_1,\tilde c_2, \tilde c_3, \tilde c_4$ are constants which clearly can be estimated above by a constant depending only on $N, m, p, a,     D_1, D_2,$ $  \delta , c_m,  \|g_1\|_{c^{m-1,1}(\overline{ W})}$, $  \|g_2\|_{c^{m-1,1}(\overline{ W})}$.
Thus, $T$ maps $W^{m,p}(W\times ]a,\infty [)$
to ${W^{m,p}({\mathcal{O}}_2)}$, and is a linear and continuous map with $\| T \|$ as in the statement.

The argument above works also for the
case $p=\infty$ provided that integrals  are replaced by the corresponding $L^{\infty}$ norms.

Finally, if $v\in L^1_{loc}(W\times ]a,\infty [)$ is such that $v=0 $ on ${\mathcal{O}}_1^c$ then $v=0 $ on ${\mathcal{O}}_{1,c_k}^c$
hence $v(\Phi_{c_k})=0 $ on ${\mathcal{O}}_2^c $ for all $k=1,\dots , m$; thus $T[v] =0$ on $ {\mathcal{O}}_2^c $.
\hfill $\Box$

\vspace{20pt}

{\bf Proof of Theorem~\ref{pretran}.} We divide the proof into three steps.

{\it Step 1.}
Recall  that $r_j(V_j)$ is a cuboid
$\Pi_{k=1}^N]a_{kj},b_{kj} [$  and
\begin{equation}
\label{sub}
r_j(\Omega_i\cap V_j)=\left\{(\bar x,x_N)\in\rea^N: \bar x\in \Pi_{k=1}^{N-1}
]a_{kj},b_{kj} [,\ a_{Nj} <x_N<g_{i,j}(\bar x) \right\},
\end{equation}
for $i=1,2$, where $g_{i,j}\in C^{m-1,1}({\overline{W}}_j)$  and
\begin{equation}a_{Nj}+
\frac{\rho}{2} <  g_{2,j}(\bar x ),\   g_{1,j}(\bar x) <
b_{Nj} -\frac{\rho}{2} ,
\end{equation}
for all $j=1,\dots ,s'$, $\bar x\in \overline{W}_j$.

For each
$j=1,\dots ,s'$ we  apply Lemma~\ref{paste},
with $W=W_j$, $a=a_{Nj}$, $D_1=b_{Nj}-\rho /2$, $D_2=a_{Nj} + \rho /2 $, $g_1=g_{1,j}$,
$g_2=g_{2,j}$, hence
\begin{equation}
\label{sub2}  \delta =\frac{1}{2} \min_{j=1, \dots , s'} \frac{\rho}{2(b_{Nj}-a_{Nj}-\rho)  }\, ,
\end{equation}
and with $c_k=k-1 +1/\delta $, for all $k=1, \dots , m$, and $\delta_k $ determined by (\ref{paste1}).
Accordingly, for each $j=1,\dots ,s'$, we consider the sets
${\mathcal{O}}_1 ={\mathcal{O}}_{1,j} $,
${\mathcal{O}}_{2} ={\mathcal{O}}_{2,j} $,
${\mathcal{O}}_{3} ={\mathcal{O}}_{3,j} $ defined by (\ref{graf2}) and the map
$T=T_{j}$, $T_j:L^1_{loc }(W_j\times ]a_{Nj},\infty [)\to L^1_{loc }( {\mathcal{O}}_{2,j} )$  defined by (\ref{paste2}). Observe that
${\mathcal{O}}_{i,j}=r_j (\Omega_i\cap V_j)$, $i=1,2$. Finally, for all
$j=1,\dots ,s'  $, we set
\begin{equation}
\label{proj}
{\mathcal{T}}_{j}[v]\equiv  (T_{j}[  (v\circ r_j^{(-1)})_{|_{W_j\times ]a_{Nj}  ,\infty [}} ]  )\circ r_j\, ,
\end{equation}
for all $v\in L^1_{loc}({\mathbb{R}}^N)$, and
\begin{equation}
\label{trej}
\Omega_{3,j}\equiv r_j^{(-1)} ({\mathcal{O}}_{3 ,j} ) .
\end{equation}

By Lemma~\ref{paste} it follows that
 ${\mathcal{T}}_j: L^1_{loc}({\mathbb{R}}^N )\to L^1_{loc}(\Omega_2\cap V_j)$,
 ${\mathcal{T}}_j: W^{m,p}  ({\mathbb{R}}^N  )\to W^{m,p} (\Omega_2\cap V_j) $,
 for all $1\le p\le \infty$, $j=1, \dots , s'$. Moreover,
 ${\mathcal{T}}_j[v](x)=v(x)$ for all $v\in L^1_{loc }({\mathbb{R}}^N)$,  $x\in \Omega_{3,j}$,
and there exists $C_{1,j}>0$ depending only on
$N,m, \rho, a_{Nj}, b_{Nj} , M$  and there exists $C_{2,j}>0$ depending only on $ \rho, a_{Nj}, b_{Nj} $ such that
$ \| {\mathcal{T}}_j \|\le C_{1,j}$ and  $|(\Omega_2\cap V_j)  \setminus   \Omega_{3,j}| \le C_{2,j} |(\Omega_1\cap V_j)\setminus  (\Omega _2\cap V_j)|$. Furthermore, if $v\in L^1_{loc}({\mathbb{R}}^N)$ and  $v=0 $ on $(\Omega_1\cap V_j)^c$ then $ {\mathcal{T}}_jv=0 $ on $(\Omega_2\cap V_j)^c$.

{\it Step 2.} We paste together the functions $ {\mathcal {T}}_j$ defined in {\it Step 1}.
To do so, we consider a partition of unity
$\{\psi_j\}_{j=1}^{s}$ such that $\psi_j\in C_{c}^{\infty}((V_j)_{\frac{3}{4}\rho})$
for all $j=1,\dots , s$ and such that $\sum_{j=1}^s\psi_j(x)=1$,
$0\le \psi_j(x)\le 1$ and $|\nabla \psi_j(x)|\le C_3$
for all
$x\in \cup_{j=1}^s(V_j)_{\rho}$,
where $C_3$ depends only on
$ {\mathcal{A}}$.

For all $j=1,\dots ,s'$, let ${\mathcal{T}}_j:L^1_{loc}({\mathbb{R}}^N )\to L^1_{loc}(\Omega_2\cap V_j)$
be as in {\it Step 1}, and for  all $s' <j\le s$, let
${\mathcal{T}}_j$ be the restriction operator  from   $L^1_{loc}({\mathbb{R}}^N )   $  to $L^1_{loc}(V_j)$.
Then we consider the operator ${\mathcal{T}} $ of $L^1_{loc}  ({\mathbb{R}}^N )$ to
$L^1_{loc}(\Omega_2)$  which takes $v\in L^1_{loc}({\mathbb{R}}^N )$ to
\begin{equation}
\label{tauuu}
{\mathcal{T}} [v] = \sum _{j=1}^s{\mathcal{T}}_j[\psi_j v],
\end{equation}
for all $v\in  L^1_{loc}({\mathbb{R}}^N ) $. Clearly, if $v\in W^{m,p}({\mathbb{R}}^N)$ then ${\mathcal{T}}_j[\psi_j v]\in W^{m,p}(\Omega_2)$.

{\it Step 3.} Since $\Omega_1\in C^{m-1,1}_M({\mathcal {A}})$, by Burenkov~\cite[Thm.~3,~p.~285]{bu} there exists a linear extension operator
$$
E_{\mathcal{N}}:W^{m,p}(\Omega_1)\to W^{m,p}({\mathbb{R}}^N)
$$
with $\| E_{\mathcal{N}} \|$ depending only on  $ {\mathcal{A}}, m$.
Let  $$E_{\mathcal{D}}: W^{m,p}_0(\Omega_1)\to W^{m,p}({\mathbb{R}}^N)$$ be the extension-by-zero operator.

We set
\begin{equation}
{\mathcal{T}}_{{\mathcal{D}}}[u]={\mathcal{T}}[E_{\mathcal{D}}u],\ \ \ \ {\mathcal{T}}_{{\mathcal{N}}}[v]={\mathcal{T}}[E_{\mathcal{N}}v]
\end{equation}
for all $u\in W^{m,p}_0(\Omega )$, $v\in W^{m,p}(\Omega )$,  and
\begin{equation}\label{omegatre}\Omega _3= \Omega _1\setminus (
\cup _{j=1}^{s'}(\Omega_1 \cap V_j)\setminus \Omega_{3, j}).
\end{equation}
Note that $\Omega_3\cap V_j\subset \Omega_{3,j}$ for all $j=1,\dots , s'$, hence $\Omega_3\subset \Omega_1\cap\Omega _2$.

By {\it Step 2.} it follows that ${\mathcal{T}}_{{\mathcal{N}}}$ maps $W^{m,p}(\Omega_1)$ to $W^{m,p}(\Omega_2)$. Moreover,
 if $u\in W^{m,p}_0(\Omega_1)$   then $E_{\mathcal{D}}u$ vanishes outside $\Omega_1$, hence ${\mathcal{T}}_j[\psi_j E_{{\mathcal{D}}}u]$ vanishes outside $\Omega_2 $  and  ${\mathcal{T}}_{{\mathcal{D}}}[u]\in W^{m,p}_0(\Omega_2)$.  Thus ${\mathcal{T}}_{{\mathcal{D}}}$ maps $W^{m,p}_0(\Omega_1)$ to $W^{m,p}_0(\Omega_2)$.

Statement {\it (i)} follows by {\it Step 1} and by the properties of the extension operators $E_{\mathcal{N}}$, $E_{\mathcal{D}}$.  The equalities in (\ref{identity}) immediately follow by Lemma~\ref{paste}. Finally,
 inequality (\ref{measure})
can be deduced by Lemma~\ref{graf} (i)  by using exactly the same argument in the proof of \cite[Lemma~4.23]{buladir}. \hfill $\Box $

\section{Sharp estimates for the variation of the eigenvalues  via the Lebesgue measure}

In this section we prove stability estimates for the eigenvalues $\lambda _{n,{\mathcal{D}}}[\Omega ] $, $\lambda _{n,{\mathcal{N}}}[\Omega ] $ defined in Definition~\ref{defset}.
Recall that $\lambda _{n,{\mathcal{D}}}[\Omega ] $, $\lambda _{n,{\mathcal{N}}}[\Omega ] $  are the eigenvalues of the operators
$H_{W^{m,2}_0(\Omega )} $,  $H_{W^{m,2}(\Omega )} $ respectively.

By $\varphi_{n,{\mathcal{D}}}[\Omega ]$ and  $ \varphi_{n,{\mathcal{N}}}[\Omega ]$ we denote a sequence of orthonormal eigenfunctions corresponding to
$\lambda _{n,{\mathcal{D}}}[\Omega ] $ and  $\lambda _{n,{\mathcal{N}}}[\Omega ] $ respectively.

When no distinction between the Dirichlet and the Neumann case is required and we refer to both, we simply write $\lambda_n[\Omega ]$, $\varphi_n[\Omega ]$, $H_{\Omega }$  to indicate the eigenvalues and the corresponding eigenfunctions and operators.

The following statement hold for both Dirichlet and Neumann boundary conditions.

\begin{thm}
\label{sharp0}
Let
${\mathcal{A}}= (  \rho , s,s', \{V_j\}_{j=1}^s, \{r_j\}_{j=1}^{s} ) $ be an atlas in ${\mathbb{R}}^N$, $m\in {\mathbb{N}}$, $ M, \theta >0$.
 For all $\alpha,\beta\in {\mathbb{N}}_0^N$ with $|\alpha |=|\beta |=m$, let $A_{\alpha\beta }$ be measurable
real-valued functions defined on $\cup_{j=1}^sV_j$, satisfying
$A_{\alpha \beta } =A_{\beta \alpha }$  and condition (\ref{elp}).

Let $2<p\le \infty $, $0< M_n <\infty$ for all  $n\in\mathbb{N}$,
and $ {\mathfrak{A}}=\left\{\Omega \in C^{m-1,1}_M({\mathcal{A}}):\ \right.$ $  \left. \|\varphi_n[\Omega ] \|_{W^{m,p}(\Omega )}\le
M_n\ for\ all\ n\in\mathbb{N} \right\}$.

Then for each $n\in\mathbb{N}$ there exists $c_n >0$
depending only on $n,  {\mathcal{A}},  m,  M, \theta , p$, $M_1, \dots ,M_n$  such
that
\begin{equation}
\label{sharp1} \lambda_{n}[\Omega_2]\le \lambda_{n}[\Omega_1]+ c_n
|\Omega_1 \vartriangle \Omega_2 |^{1-\frac{2}{p}},
\end{equation}
for all $\Omega_1\in {\mathfrak{A}}$, $\Omega_2\in C^{m-1,1}_M({\mathcal{A}} )$ such that
$|\Omega_1 \vartriangle \Omega_2 |< c_n^{-1} $ .

\end{thm}

{\bf Proof.}
Let $\Omega_1  \in {\mathfrak{A}}$ and $\Omega_2 \in C^{m-1,1}_M({\mathcal{A}} )$. To shorten our notation we set $\varphi_{n,1}=\varphi_n[\Omega_1]$, for all $n\in {\mathbb{N}}$. We denote by ${\mathcal{L}}_1$ the space of the finite linear combinations
of the eigenfunctions $\varphi_{n,1}$.
Moreover, we  define a linear operator $$T_{12}:{\mathcal{L}}_1\to {\rm {Dom}} (H^{1/2}_{\Omega_2})$$ by setting in the Dirichlet case
$$
T_{12}[\varphi_{n,1}]={\mathcal{T}}_{{\mathcal{D}}}\varphi_{n,1},
$$
and in the Neumann case
$$
T_{12}[\varphi_{n,1}]={\mathcal{T}}_{{\mathcal{N}}}\varphi_{n,1}.
$$
for all $n\in {\mathbb{N}}$.
Here
$$
{\mathcal{T}}_{D}:W^{m,p}_0(\Omega_1)\to W^{m,p}_0(\Omega_2)\ \ {\rm and}\ \ {\mathcal{T}}_{N}:W^{m,p}(\Omega_1)\to W^{m,p}(\Omega_2)
$$
are the operators provided by Theorem~\ref{pretran}. Note that $T_{12}$ is well-defined. Indeed, by assumption ${\mathcal{L}}_1 \subset W^{m,p}(\Omega_1)$, and in the Dirichlet case ${\mathcal{L}}_1\subset W^{m,p}_0(\Omega_1)$. Moreover,
$T_{12}$ takes values in ${\rm {Dom}} (H^{1/2}_{\Omega_2})$ because in the  Dirichlet case $W^{m,p}_0(\Omega _2)$ $ \subset$ $ W^{m,2}_0(\Omega _2)={\rm {Dom}} (H^{1/2}_{\Omega_2})$, and  in the Neumann case
   $W^{m,p}(\Omega _2)$ $\subset W^{m,2}(\Omega _2)={\rm {Dom}} (H^{1/2}_{\Omega_2})$.

To prove (\ref{sharp1}) we apply the general spectral stability theorem \cite[Thm.~3.2]{bula}. In the terminology of \cite{bula}, we need to prove that $T_{12}$ is a `transition operator' from $H_{\Omega_1}$ to $H_{\Omega_2}$. To do so, we prove  inequalities (\ref{tran1}) and (\ref{tran2}) below.

By  Theorem \ref{pretran}, $T_{12}\varphi_n=\varphi_n$ on $\Omega_3$ where $\Omega_3$ is as in Theorem~\ref{pretran} $(ii)$. Thus
\begin{eqnarray}\label{addit}
\label{sharp4} \lefteqn{
(H_{\Omega_2}^{1/2}T_{12}\varphi_{k,1},H_{\Omega_2}^{1/2}T_{12}\varphi_{l,1})_{L^2(\Omega_2)}
=Q_{\Omega_2}(T_{12}\varphi_{k,1},T_{12}\varphi_{l,1})
}\\
& & =
Q_{\Omega_3}( T_{12}\varphi_{k,1},T_{12}\varphi_{l,1} )
+ Q_{\Omega_2\setminus
\Omega_3} ( T_{12}\varphi_{k,1},T_{12}\varphi_{l,1} )
\nonumber\\ & &
=
Q_{\Omega_3}( \varphi_{k,1}, \varphi_{l,1} )
+ Q_{\Omega_2\setminus
\Omega_3} ( T_{12}\varphi_{k,1},T_{12}\varphi_{l,1} )
\nonumber\\
& & =
Q_{\Omega_1}( \varphi_{k,1}, \varphi_{l,1} )-
Q_{\Omega_1\setminus \Omega _3}( \varphi_{k,1}, \varphi_{l,1} )
+Q_{\Omega_2\setminus
\Omega_3} ( T_{12}\varphi_{k,1},T_{12}\varphi_{l,1} )
\nonumber\\
& & =
(H_{\Omega_1}^{1/2}\varphi_{k,1},H_{\Omega_1}^{1/2}\varphi_{l,1})_{L^2(\Omega_1)}
-
Q_{\Omega_1\setminus \Omega _3}( \varphi_{k,1}, \varphi_{l,1} )
+Q_{\Omega_2\setminus
\Omega_3} ( T_{12}\varphi_{k,1},T_{12}\varphi_{l,1} ),\nonumber
\end{eqnarray}
for all $k,l\in\mathbb{N}$. By H\"{o}lder's inequality
\begin{eqnarray}
Q_{\Omega_1\setminus \Omega _3}( \varphi_{k,1}, \varphi_{l,1} )
\le c  M_{k} M_{l}
|\Omega_1\setminus\Omega_3|^{1-\frac{2}{p}}
\end{eqnarray}
and by Theorem~\ref{pretran} we have

\begin{eqnarray}
Q_{\Omega_2\setminus
\Omega_3} ( T_{12}\varphi_{k,1},T_{12}\varphi_{l,1} )
 \le
c  M_{k}M_{l}|\Omega_2\setminus
\Omega_3|^{1-\frac{2}{p}},
\end{eqnarray}
and
\begin{equation}
\label{sharp2,5} |\Omega_1\setminus \Omega_3|,\ |\Omega_2\setminus \Omega_3|\le c|\Omega_1\vartriangle \Omega_2| ,
\end{equation}
where $c>0$ depends only on $ {\mathcal{A}}, m, M, \theta , p  $.
Thus by (\ref{sharp4})-(\ref{sharp2,5}) it follows that
\begin{eqnarray}
\label{tran1,5}
\lefteqn{
|(H_{\Omega_2}^{1/2}T_{12}\varphi_{k,1},H_{\Omega_2}^{1/2}T_{12}\varphi_{l,1}  )_{L^2(\Omega_2 )}} \nonumber  \\ & & \qquad\qquad - (H_{\Omega_1}^{1/2}\varphi_{k,1},H_{\Omega_1}^{1/2}\varphi_{l,1}  )_{L^2(\Omega_1 )}   |  \le  \tilde c_5M_kM_l
|\Omega_1\vartriangle \Omega_2|^{1-\frac{2}{p}},  \label{tran1}
\end{eqnarray}
and similarly
\begin{equation}\label{tran2}
|(T_{12}\varphi_{k,1},T_{12}\varphi_{l,1}  )_{L^2(\Omega_2 )}
 - (\varphi_{k,1},\varphi_{l,1}  )_{L^2(\Omega_1 )}   |
\le \tilde c_6 M_kM_l
|\Omega_1\vartriangle \Omega_2|^{1-\frac{2}{p}},
\end{equation}
for all $k,l\in {\mathbb{N}}$, where $\tilde c_5,\tilde c_6>0$ depend only on $ {\mathcal{A}}, m, M, \theta , p  $.

By (\ref{tran1,5}), (\ref{tran2}) it follows that $T_{12}$ is a transition operator from $H_{\Omega_1}$ to $H_{\Omega_2}$
with parameters $a_{kl}=\tilde c_5M_kM_l$, $b_{kl}=\tilde c_6M_kM_l$ and measure of vicinity $\delta (H_{\Omega_1},H_{\Omega_2})=|\Omega_1\vartriangle \Omega_2|^{1-\frac{2}{p}}$ (see \cite[Def.~3.1]{bula}). Thus by \cite[Thm.~3.2]{bula} it follows that
\begin{equation}
\label{burlamest}
\lambda_n[\Omega _2]\le \lambda_n[\Omega_1]+(2a_n\lambda_n[\Omega_1]+b_n)\delta (H_{\Omega_1},H_{\Omega_2}).
\end{equation}
if $\delta (H_{\Omega_1},H_{\Omega_2})\le (2a_n)^{-1}$, where $a_n=(\sum_{k,l=1}^na_{kl}^2)^{1/2}=\tilde c_5\sum_{k=1}^nM_k^2$, $b_n=(\sum_{k,l=1}^nb_{kl}^2)^{1/2}=\tilde c_6\sum_{k=1}^nM_k^2$.
Furthermore, by \cite[Lemma~3.2]{bulahigh} there exists $\Lambda_n>0$ depending only on $n, {\mathcal{A}}, m, \theta$ such that
\begin{equation}\label{unifest}
\lambda_n[\Omega ]\le \Lambda_n
\end{equation}
   for all $\Omega \in C^{m-1,1}_{M}({\mathcal{A}})$. Thus, inequality (\ref{sharp1}) follows by combining (\ref{burlamest}) and (\ref{unifest}).
\hfill $\Box$\\

\begin{rem}
It can be traced that starting with (\ref{addit}) one can obtain the estimate
$$
\lambda_n[\Omega _2]\le \lambda_n[\Omega _1]+c_n \sum_{k=1}^n\| \varphi_k [\Omega_1]\|_{W^{m,2}(\Omega_1\vartriangle \Omega_2)}
$$
which in some cases (depending on the properties of $\varphi_1[\Omega_1],\dots ,\varphi_n[\Omega_1]$ near the boundary of $\Omega_1$) can be better than estimate (\ref{sharp1}).
\end{rem}

It is well-known that if $\Omega_2\subset \Omega_1$ then
$\lambda_{n,{\mathcal{D}}}[\Omega_1]\le \lambda_{n,{\mathcal{D}}}[\Omega_2]$. Thus by
Theorem~\ref{sharp0} we immediately deduce the following corollary concerning Dirichlet eigenvalues (for the proof of the sharpness of estimate (\ref{sharpincest}), see Section 5).

\begin{corol}
\label{sharpinc}
Let
${\mathcal{A}}= (  \rho , s,s', \{V_j\}_{j=1}^s, \{r_j\}_{j=1}^{s} ) $ be an atlas in ${\mathbb{R}}^N$, $m\in {\mathbb{N}}$, $ M, \theta >0$.
Let $\Omega_1\in C^{m-1,1}_{M}({\mathcal{A}})$.
 For all $\alpha,\beta\in {\mathbb{N}}_0^N$ with $|\alpha |=|\beta |=m$, let $A_{\alpha\beta }$ be measurable
real-valued functions defined on $\Omega_1$, satisfying
$A_{\alpha \beta } =A_{\beta \alpha }$  and condition (\ref{elp}).

Assume that $2<p\le\infty $ and $ \varphi_{n ,{\mathcal{D}}}[\Omega_1]\in W^{m,p}(\Omega_1)$ for all
$n\in\mathbb{N}$. Then for each $n\in\mathbb{N}$ there exists
$c_{n}>0$ depending only on $n, {\mathcal{A}}, m, M,$ $ \theta, p,$ $  \|\varphi_k[\Omega_1] \|_{W^{m,p}(\Omega_1)}$ $k=1,\dots ,n$,
 such that
 \begin{equation}
 \label{sharpincest}
 \lambda_{n,{\mathcal{D}}}[\Omega_1]\le \lambda_{n,{\mathcal{D}}}[\Omega_2]\le \lambda_{n,{\mathcal{D}}}[\Omega_1]+c_n|\Omega_1\setminus \Omega_2|^{1-\frac{2}{p}},
 \end{equation}
for all $\Omega_2$ of class $ C^{m-1,1}_M({\mathcal{A}} )$ such that $\Omega_2\subset
\Omega_1$ and $|\Omega_1\setminus\Omega_2|< c_n^{-1}$.

Moreover, in general the exponent $1-\frac{2}{p}$ in (\ref{sharpincest}) cannot be replaced by $1-\frac{2}{p}+\delta$ where $\delta >0$ is a constant independent of $p$.
\end{corol}

If we assume that
both $\Omega_1$ and $\Omega_2$ belong to ${\mathfrak{A }}$ then it is possible to swap $\Omega_1$ and $\Omega_2$ in (\ref{sharp1}). In this way we obtain a two-sided estimate for both Dirichlet and Neumann eigenvalues without assuming that $\Omega_2\subset \Omega_1$ as in Corollary \ref{sharpinc}.

\begin{corol}
\label{unifbounded}
Let
${\mathcal{A}}= (  \rho , s,s', \{V_j\}_{j=1}^s, \{r_j\}_{j=1}^{s} ) $ be an atlas in ${\mathbb{R}}^N$, $m\in {\mathbb{N}}$, $ M, \theta >0$.
 For all $\alpha,\beta\in {\mathbb{N}}_0^N$ with $|\alpha |=|\beta |=m$, let $A_{\alpha\beta }$ be measurable
real-valued functions defined on $\cup_{j=1}^sV_j$, satisfying
$A_{\alpha \beta } =A_{\beta \alpha }$  and condition (\ref{elp}).

Let $2<p\le\infty$ and  let ${\mathfrak{A}}$ be a family of open sets of
class $C^{m-1,1}_M({\mathcal{A}})$ such that for each $n\in\mathbb{N}$
condition (\ref{unifbound}) is satisfied.

Then for each $n\in\mathbb{N}$ there exists $c_n>0$ depending only on  $n, {\mathcal{A}}, m, M, \theta,$ $  p$, $\sup_{\Omega \in {\mathfrak{A}} }\|\varphi_k[\Omega ] \|_{W^{m,p}(\Omega)}$ $k=1,\dots ,n$,   such
that
\begin{equation}
\label{unifboundedest}
|\lambda_n[\Omega_1]-\lambda_n[\Omega_2] |\le c_n |\Omega_1
\vartriangle \Omega_2 |^{1-\frac{2}{p}},
\end{equation}
for all $\Omega_1,\Omega_2\in {\mathfrak{A}}$ such that $ |\Omega_1
\vartriangle \Omega_2 |< c_n^{-1}$.

\end{corol}

If ${\mathfrak{A}}$ is a family of open sets with sufficiently smooth boundaries then condition (\ref{unifbound}) is satisfied with $p=\infty $.

\begin{lem}
\label{smoothestlem}
Let
${\mathcal{A}}= (  \rho , s,s', \{V_j\}_{j=1}^s, \{r_j\}_{j=1}^{s} ) $ be an atlas in ${\mathbb{R}}^N$, $m\in {\mathbb{N}}$, $B, M, \theta >0$.
 For all $\alpha,\beta\in {\mathbb{N}}_0^N$ with $|\alpha |=|\beta |=m$, let $A_{\alpha\beta }\in C^{m}(\overline{\cup_{j=1}^sV_j })  $ satisfy
$A_{\alpha \beta } =A_{\beta \alpha }$, $\|A_{\alpha \beta } \|_{c^{m}(\overline{\cup_{j=1}^sV_j })}\le B  $,  and condition (\ref{elp}).
Then $\varphi_n[\Omega ]\in W^{2m-1,\infty }(\Omega ) $  and  there exists $C>0$ depending only on ${\mathcal{A}}, m, B, M, \theta$ such that
\begin{equation}
\label{inftyeiglem}
\| \varphi_n[\Omega ]\|_{W^{k,\infty }(\Omega )}\le C (1+\lambda_n[\Omega ])^{\frac{N}{4m}+\frac{k}{2m}}
\end{equation}
for all $k=0,\dots , 2m-1$ and  $\Omega\in  C^{2m}_M({{\mathcal{A}}}) $.
\end{lem}

{\bf Proof.}
It is well-known that under our regularity assumptions ${\rm Dom }(H)\subset W^{2m,2}(\Omega )$ (see {\it e.g.}, Agmon~\cite[Sec.~9]{agmonlect}).
Moreover, since the coefficients $A_{\alpha \beta }$ are of class $C^m$ and we impose either Dirichlet or Neumann boundary conditions, we can resort to the general setting of Agmon~\cite{agmon62} (see \cite[pp.~141-143]{agmonlect} for details).

Thus, by \cite[Thm.~1.1~and~the~Lemma on p.131]{agmon62} it follows that if $u\in {\rm Dom }(H)$  and $Hu\in L^p(\Omega )$ for some $p>1$ then
$u\in W^{2m, p}(\Omega )$ and
\begin{equation}\label{aprioriest}
\|u\|_{W^{2m,p}(\Omega )}\le c( \|Hu\|_{L^p(\Omega )}+\|u\|_{L^{p}(\Omega )} ),
\end{equation}
where $c$ is a positive constant.
In particular if $\varphi $ is an eigenfunction corresponding to an eigenvalue $\lambda$ and $\varphi \in L^p(\Omega )$ then
\begin{equation}\label{aprioriesteige}
\|u\|_{W^{2m,p}(\Omega )}\le c(1+\lambda )\|u\|_{L^p(\Omega )}.
\end{equation}

By the apriori estimate (\ref{aprioriesteige}) and a  bootstrap argument one can finally prove estimate (\ref{inftyeiglem}) (see for instance
Burenkov and Lamberti~\cite[Thm.~5.1]{bula} where in the proof one has simply to replace \cite[(5.5)]{bula} by (\ref{aprioriesteige})).
\hfill $\Box$\\

By Corollary~\ref{unifbounded} and Lemma~\ref{smoothestlem} we immediately deduce the validity of the following

\begin{corol}
\label{smoothest}
Let
${\mathcal{A}}= (  \rho , s,s', \{V_j\}_{j=1}^s, \{r_j\}_{j=1}^{s} ) $ be an atlas in ${\mathbb{R}}^N$, $m\in {\mathbb{N}}$, $B, M, \theta >0$.
 For all $\alpha,\beta\in {\mathbb{N}}_0^N$ with $|\alpha |=|\beta |=m$, let $A_{\alpha\beta }\in C^{m} (\overline{\cup_{j=1}^sV_j })  $ satisfy
$A_{\alpha \beta } =A_{\beta \alpha }$, $\|A_{\alpha \beta } \|_{c^{m}(\overline{\cup_{j=1}^sV_j })}\le B  $,  and condition (\ref{elp}).
Then for all $n\in {\mathbb{N}}$  there exists $c_n>0$ depending only on $n, {\mathcal{A}}, m,B,  M, \theta$ such that
\begin{equation}
\label{inftyeig}
|\lambda_n[\Omega_1]-\lambda_n[\Omega_2]|\le c_n |\Omega_1\vartriangle \Omega_2  | ,
\end{equation}
for all  $\Omega_1,\Omega_2\in  C^{2m}_M({{\mathcal{A}}}) $ satisfying $|\Omega_1\vartriangle \Omega_2|<c_n^{-1}$.
\end{corol}

\section{An example}

We consider an example which proves that in the class of Lipschitz domains the exponent in estimates (\ref{intro1}) and (\ref{intro2}) cannot, in general,  be larger than $1-2/p$. For this purpose we consider the Dirichlet and Neumann Laplacians  on the circular sector  $\Omega \subset {\mathbb{R}}^2$ of radius $R=1$ and angle $2\beta $ with $0<\beta <\pi$.
In polar coordinates
\begin{equation}
\label{omega}
\Omega =\{(\rho , \theta ):\ 0<\rho <1, \ -\beta <\theta < \beta   \}.
\end{equation}
For $0<\epsilon  <1$ we consider the deformation $\Omega (\epsilon )$ of $\Omega $ given by
\begin{equation}
\label{omegaepsilon}
\Omega (\epsilon )=\{(\rho , \theta ):\ \epsilon <\rho <1, \ -\beta <\theta < \beta   \}.
\end{equation}
Here we are interested in the behavior of the eigenvalues   of the Dirichlet and Neumann Laplacians on $\Omega (\epsilon ) $
as $\epsilon \to 0$.

In the case of the Dirichlet Laplacian on $\Omega$ all the eigenvalues are the positive solutions of the equations

\begin{equation}
\label{zerojei}
J_{\nu}( \sqrt {\lambda })=0,
\end{equation}
where $J_{\nu }$ is the Bessel function of the first kind and order $\nu$, with $\nu = \pi k /(2\beta )$, $k\in {\mathbb{N}}$.

Note that $\nu >1/2$ for all $0<\beta <\pi $, $k\in {\mathbb{N}}$, and that $\nu <1  $ if an only if $k=1$ and $\pi /2<\beta <\pi$.

For our purposes, it is enough to restrict our attention to the case  $\nu \notin {\mathbb{N}}$: in this case the eigenvalues of the Dirichlet Laplacian on $\Omega (\epsilon )$  are the positive solutions of the cross-product equations

\begin{equation}
\label{crossproduct}
J_{\nu }(\sqrt{\lambda })J_{-\nu }(\epsilon \sqrt{\lambda }))-J_{-\nu }(\sqrt{\lambda })J_{\nu }(\epsilon \sqrt{\lambda })=0.
\end{equation}

Recall that for a Bessel function of the first kind and order $\mu $ ($\mu\ne -1, -2, \dots  $) we have $J_{\mu }(s)= s^{\mu }H_{\mu }(s^2)$, $s\in {\mathbb{R}}$, where $H_{\mu }$ is an analytic function such that $H_{\mu }(0)\ne 0$, see (\cite[\S 9.1.10]{olv}).

Assume that $\lambda_{*}$ is a fixed eigenvalue of the Dirichlet  Laplacian on $\Omega $, {\it i.e.,} $\lambda_{*}$ is a fixed zero of
$H_{\nu }$. It is known that $H_{-\nu }(\lambda _{*})\ne 0$.
Thus, in a sufficiently small small neighborhood of $\lambda_{*}$ and for sufficiently small $\epsilon \geq 0$, equation (\ref{crossproduct})
can be rewritten as
\begin{equation}
\label{effeeq}
f(\lambda  )-\epsilon ^{2\nu } f(\epsilon^2\lambda  )=0,
\end{equation}
where $f(\lambda )=H_{\nu }(\lambda )/ H_{-\nu }(\lambda )$ is an analytic function in a neighborhood of zero and in a neighborhood of
$\lambda_*$.

It is immediate to verify that if $\epsilon=0$ then the positive solutions of (\ref{zerojei}) coincide with the positive solutions of equation (\ref{effeeq}).
Thus, {\it for each $0\le \epsilon <1$
  the eigenvalues $\lambda $  of the Dirichlet Laplacian  on $\Omega (\epsilon )$ are exactly the zeros of equation (\ref{effeeq})}  (here it is understood that $\Omega (0)=\Omega$).

We set $\delta =\epsilon^{2\nu }  $, so that equation (\ref{effeeq}) can be rewritten as
\begin{equation}\label{effedelta}
f(\lambda )-\delta f(\delta ^{\frac{1}{\nu }}\lambda )=0.
\end{equation}

 Observe that the left-hand side of equation (\ref{effedelta}) defines   a function of class $C^{1}$ in the variables $\delta , \lambda $, for all $(\delta , \lambda )$ a neighborhood of the point  $(0, \lambda _*)$.  Note that
 $H'_{\nu }(\lambda_*)\ne 0$ since all positive zeros of the  Bessel functions $J_{\nu }$ are simple, see \cite[9.5.2]{olv}. Thus, $f'(\lambda_*)\ne 0$ and  and by the Implicit Function Theorem the zeros of equation (\ref{effedelta}) in a neighborhood of $(0,\lambda_*)$ are given by the
graph of a function $\delta \mapsto \lambda (\delta  )$ of class $C^1 $ such that $\lambda (0)=\lambda_*$. Moreover, since the derivative
of  $f(\lambda )-\delta f(\delta ^{\frac{1}{\nu }}\lambda )$ with respect to  $\delta $ at the point $(0,\lambda_*)$ is equal to $-f(0)$ then by the Implicit Function Theorem we have that
\begin{equation}\label{expans}
\lambda (\delta )= \lambda_*+\frac{f(0)}{f'(\lambda_*)}\delta +o(\delta ), \ \ {\rm as}\ \ \delta  \to 0^+.
\end{equation}
Note that  $f(0)\ne 0$. This clearly implies that
\begin{equation}\label{asymp}
|\lambda (\epsilon )-\lambda _* |=C |\Omega \setminus \Omega (\epsilon )|^{\nu }+ o(|\Omega \setminus \Omega (\epsilon )|^{\nu }),\ \ \ {\rm as}\ \ \epsilon \to 0^+,
\end{equation}
where $C$ is a positive constant.

We note that the eigenspace  of the Dirichlet Laplacian on $\Omega $ corresponding to the eigenvalue $\lambda _*$ is spanned by the function $U$
defined in polar coordinates by
\begin{equation}\label{expliciteigen}
U(\rho, \theta )= J_{\nu }(\rho \sqrt {\lambda_* })\sin \nu (\theta +\beta ) = (\rho\sqrt {\lambda _*})^{\nu }H_{\nu }(\rho^2\lambda_*)\sin  \nu (\theta +\beta ) ,
\end{equation}
for all $0<\rho < 1$, $-\beta <\theta <\beta $.
Clearly, $U\in L^{\infty }(\Omega )$ as expected, and  if $\nu \geq 1  $ then $\nabla U\in L^{\infty }(\Omega ) $, whilst if $1/2<\nu <1$ then
$\nabla U\in L^{p}(\Omega )$ if and only if  $1\le p<2/(1-\nu ) $.
Thus by applying  estimate  (\ref{intro2}) we obtain that
if $1/2 < \nu <1$ then  for any $ 0<\gamma <\nu$ there exists $c_{\gamma }>0$ such that
\begin{equation}\label{pest}
| \lambda (\epsilon )-\lambda_* |\le c_{\gamma } |\Omega \setminus \Omega (\epsilon )|^{\gamma },
\end{equation}
if $ |\Omega \setminus \Omega (\epsilon )|<c_{\gamma }^{-1}$, whilst
if $\nu \geq  1$ then there exists $c>0$ such that
\begin{equation}\label{inftyest}
| \lambda (\epsilon )-\lambda_* |\le c |\Omega \setminus \Omega (\epsilon )|,
\end{equation}
if $ |\Omega \setminus \Omega (\epsilon )|<c^{-1}$.

In the case of the Neumann Laplacian, equations (\ref{zerojei}), (\ref{crossproduct}) should be replaced by equations
\begin{equation}
\label{zerojeineu}
J'_{\nu }(\sqrt{ \lambda }) =0,
\end{equation}
and

\begin{equation}
\label{crossproductneu}
J_{\nu }'(\sqrt{\lambda })J_{-\nu }'(\epsilon \sqrt{\lambda }))-J_{-\nu }'(\sqrt{\lambda })J_{\nu }'(\epsilon \sqrt{\lambda })=0,
\end{equation}
 respectively. By writing $J_{\mu }'(s)=s^{\mu -1}K_{\mu }(s^2) $ where $K_{\mu }$ is a suitable analytic function not vanishing at zero, one can easily see that in the case of the Neumann Laplacian in equation (\ref{effeeq}) one should simply  replace the function $f$ by the function $g(\lambda )=K_{\nu }(\lambda )/ K_{-\nu }(\lambda )$.  Note that  $K'_{\nu }(\lambda_*)\ne 0$ since all positive zeros of  $J_{\nu }'$ are simple (use directly the Bessel equation of order $\nu$ and  Watson \cite[\S 15.3, (3), p.~486]{wat}). Thus, one can apply the same argument used above and prove
 that  (\ref{asymp}) holds also for the eigenvalues of the Neumann Laplacian.

 Note that the eigenspace  of the Neumann Laplacian on $\Omega $ corresponding to a positive  eigenvalue $\lambda _*$ is spanned by the function $V$ defined in polar coordinates by
$$
V(\rho, \theta )= J_{\nu }(\rho \sqrt {\lambda_* })\cos  \nu (\theta +\beta ) = (\rho\sqrt {\lambda _*})^{\nu }H_{\nu }(\rho^2\lambda_*)\cos  \nu (\theta +\beta ) ,
$$
for all $0<\rho <1$, $-\beta <\theta <\beta $.
Thus also for the Neumann Laplacian, we conclude that inequality (\ref{intro1}) implies  (\ref{pest}) and (\ref{inftyest}).

{\it Clearly, in both the cases of Dirichlet and Neumann boundary conditions, (\ref{asymp}) shows that if $k=1$ and $\pi /2<\beta <\pi $ ($\iff$  $1/2<\nu <1 $) then  the exponent $\gamma$ in (\ref{pest}) cannot be larger than $\nu$.
Thus, in general, the exponent in the right hand-side of estimates (\ref{intro1}), (\ref{intro2}) in the class of Lipschitz domains cannot be larger than $1-2/p$.  However,   (\ref{asymp}) and (\ref{inftyest}) also  show that for special domains and special values of the indices $n$ one may find better exponents in the right hand-side of estimates (\ref{intro1}), (\ref{intro2}).}

Note that in this example the domains $\Omega $ and $\Omega (\epsilon )$ are of class $C^{0,1}$ but not of class $C^{0,1}_M({\mathcal{A}})$
for fixed atlas ${\mathcal{A}}$ and $M>0$. In the proof below the domains $\Omega (\epsilon )$ will be modified in an appropriate way in order to define suitable domains $\tilde \Omega (\epsilon )$ belonging to the same class $C^{0,1}_M({\mathcal{A}})$.\\

{\bf Proof of the sharpness of the exponent $1-2/p$ in  (\ref{sharpincest}) for $N=2$, $m=1$, $n=1$.} In this proof, by $\lambda_{1,{\mathcal{D}}}[{\mathcal{U}}]$ we denote the first eigenvalue of the Dirichlet Laplacian  defined on a bounded domain ${\mathcal{U}}$ in ${\mathbb{R}}^2$.

Let $\Omega $ be the domain defined by (\ref{omega}) with $\pi /2<\beta <\pi$.  For all $\epsilon \in ]0,1[$ we set
$$
\tilde \Omega (\epsilon )= \left\{x=(x_1,x_2 )\in {\mathbb{R}}^2:\ g(x_2 )<x_1,\ \ |x |<1  \right\},
$$
where  $g(x_2)= \epsilon - | x_2|\tan \frac{\beta }{2}$ if $|x_2|\le \epsilon \sin \beta $, and $g( x_2)=| x_2| \cot \beta $ if
$| x_2|> \epsilon \sin \beta $.

It is easy to see that
\begin{equation}\label{inclus}
\Omega (\epsilon )\subset  \tilde \Omega (\epsilon )\subset \Omega (A\epsilon )\subset \Omega ,
\end{equation}
for all $\epsilon \in ]0,1[$,  where $\Omega (\epsilon )$ is defined by (\ref{omegaepsilon})  and   $A=\cos \frac{\beta }{2}$. By monotonicity it follows that
\begin{equation}\label{incluseigen}\lambda_{1, {\mathcal{D}}}[\Omega ]\le
\lambda_{1, {\mathcal{D}}}[\Omega (A\epsilon )]\le \lambda_{1,{\mathcal{D}}}[\tilde\Omega (\epsilon )] \le \lambda_{1, {\mathcal{D}}}[\Omega (\epsilon )].
\end{equation}
Since the eigenfunctions corresponding to the first eigenvalue of the Dirichlet Laplacian are the only eigenfunctions which do not change sign, it follows that the eigenspace corresponding to the eigenvalue $\lambda_{1,{\mathcal{D}}}[\Omega ]$ is spanned by (\ref{expliciteigen}) with $\nu =\pi / (2\beta )$.
Thus, the asymptotic behavior of $\lambda_{1,{\mathcal{D}}}[\Omega (\epsilon )]$ is given by (\ref{asymp}) with $\nu =\pi / (2\beta )$, hence
\begin{equation}
\label{main2}
\lambda _{1,{\mathcal{D}}} [\Omega (\epsilon )]= \lambda_{1, {\mathcal{D}}}[\Omega ]+ C|\Omega \setminus \Omega (\epsilon )|^{\frac{\pi}{2\beta }}+o\biggl(|\Omega \setminus \Omega (\epsilon )|^{  \frac{\pi }{2\beta}   }  \biggr),\ \ {\rm as}\ \epsilon \to 0^+.
\end{equation}
By combining (\ref{incluseigen}) and (\ref{main2}) it follows that
\begin{equation}
\label{esttilde}
C_1|\Omega \setminus \tilde \Omega (\epsilon ) | ^{\frac{\pi}{2\beta }  } \le
|\lambda_{1, {\mathcal{D}}}[\tilde \Omega (\epsilon )]-\lambda_{1, {\mathcal{D}}}[ \Omega ]|\le C_2|\Omega \setminus \tilde \Omega (\epsilon ) |  ^{\frac{\pi}{2\beta}},
\end{equation}
for all sufficiently small $\epsilon $, where $C_1,C_2$ are positive constants independent of $\epsilon$.

We now apply Corollary~\ref{sharpinc} to the Dirichlet Laplacian with $\Omega_1=\Omega $ and $\Omega _2=\tilde \Omega (\epsilon)$.
It is clear that there exists an atlas ${\mathcal{A}}$ and $M>0$ such that $\Omega$ and
$
\tilde \Omega (\epsilon )$  are of class $ C^{0,1}_M({\mathcal {A}})$
for all $\epsilon \in ]0, 1/2 [$.
Moreover, by formula (\ref{expliciteigen}) it follows that if $1\le p<4\beta /(2\beta -\pi ) $ then the eigenfunctions
$ \varphi_{n ,{\mathcal{D}}}[\Omega_1]$ of the Dirichlet Laplacian in $\Omega_1$ belong to $W^{1,p}(\Omega_1)$ for all
$n\in\mathbb{N}$. Thus the assumptions of Corollary~\ref{sharpinc}  are satisfied for such range of $p$.
Assume now by contradiction that  under the assumption of Corollary~\ref{sharpinc} estimate (\ref{sharpincest}) is valid with $|\Omega_1\setminus \Omega_2|^{1-2/p +\delta }$ replacing $|\Omega_1\setminus \Omega_2|^{1-2/p  }$, where  $\delta $ is a positive constant independent of $p$.
Since $\lim_{p\to 4\beta /(2\beta -\pi )}1-2/p=\pi /(2\beta )  $, by choosing $p$ sufficiently close to $4\beta /(2\beta -\pi ) $
it follows that the second inequality in (\ref{esttilde}) holds with $|\Omega \setminus \tilde \Omega (\epsilon ) |  ^{\frac{\pi}{2\beta}+\mu }$ replacing $|\Omega \setminus \tilde \Omega (\epsilon ) |  ^{\frac{\pi}{2\beta}}$ for some $\mu >0$
and this contradicts the first inequality in (\ref{esttilde}) as $\epsilon \to 0^+$.
\hfill $\Box$\\

{\bf Acknowledgments}: The authors are thankful to Professor M. Marletta who carried out some numerical calculations confirming the sharpness
of the exponent $1-2/p$ in estimates  (\ref{sharp1}), (\ref{sharpincest}), (\ref{unifboundedest}) and thus encouraged us to give analytic proof of the sharpness of this exponent.

We also note that further numerical calculations and some analytic computations in the spirit of Section 5 were carried out by the student N.A. Oliver under the supervision of Professor M. Marletta.

This research was supported by the research project ``Problemi di stabilit\`{a}  per operatori differenziali" of the University of Padova, Italy and by the research project PRIN 2008 ``Aspetti geometrici delle equazioni alle derivate parziali e questioni connesse''. The first author was also supported by the grant of RFBR - Russian Foundation for Basic Research (projects 09-01-00093-A, 11-01-00744-A).

The authors are also especially thankful to the anonymous Referee for critical analysis of the manuscript and valuable comments which helped to improve the presentation of the results.

\noindent

{\small
V.I. Burenkov, Faculty of Mechanics and Mathematics, L.N. Gumlyov Eurasian National University, 5 Munaitpasov Str., 010008 Astana, Kazakhstan. \\

P.D. Lamberti,  Dipartimento di Matematica Pura ed Applicata, Universit\`{a} degli Studi di Padova, Via Trieste 63, 35121 Padova, Italy.\\}

\end{document}